\documentclass[10pt,A4]{article}
\usepackage{latexsym,enumerate,amsmath,amsthm,amsopn,amstext,amscd,amsfonts,amssymb,mathrsfs,amsbsy,pictex}

\usepackage{graphicx}

\newtheorem{theorem}{Theorem}

\newtheorem{proposition}[theorem]{Proposition}

\newtheorem{rem}[theorem]{Remark}

\newtheorem{example}[theorem]{Example}

\newtheorem{lemma}[theorem]{Lemma}

\let\epsilon=\varepsilon

\let\ep=\varepsilon

\let\phi=\varphi

 \let\th=\theta

\let\tilde=\widetilde

\def\Var{{\rm Var}}

\def\tref#1{(\ref{#1})}
\newcommand{\bee}{\begin{enumerate}}
\newcommand{\eee}{\end{enumerate}}

\newcommand{\field}[1]{\mathbb{#1}}
\newcommand{\R}{\field{R}}
\newcommand{\bR}{\field{R}}
\newcommand{\bS}{\field{S}}

\renewcommand{\P}{{\rm P}}

\newcommand{\Z}{\mathscr{Z}}

\newcommand{\EE}{\field{E}}
\newcommand{\bP}{\field{P}}
\newcommand{\bE}{\field{E}}

\newcommand{\bL}{L}

\newcommand{\beqn}{\begin{equation}}
\newcommand{\eeqn}{\end{equation}}

\newcommand{\cH}{{\mathscr H}}

\def\E{{\field E}}

\def\ccl#1{{\mathscr #1}}

\def\bjk{\beta_{j\eta}}

\def\hbjk{\hat\beta_{j\eta}}

\def\psijk{\psi_{j\eta}}
\def\hf{\hat f}

\title{\sc Adaptive density estimation for directional data using needlets}
\author{\bf P.Baldi$^{\ast}$, G.Kerkyacharian$^{\ast\ast}$,\\\bf  D.Marinucci$^{\ast}$ and
D.Picard$^{\ast\ast}$\\
\\
$^{\ast}$Dipartimento di Matematica,
Universit\`a di Roma {\sl Tor Vergata},
Italy\\
$^{\ast\ast}$Laboratoire de Probabilit\'es et Mod\`eles
Al\'eatoires, Paris, France}
\date{}
\begin{document}
\maketitle
\begin{abstract}
This paper is concerned with density estimation of directional data on the sphere.
We introduce a procedure based on thresholding on a new type of
spherical wavelets called {\it needlets}. We establish a minimax result and prove
its optimality. We are motivated by astrophysical applications,
in particular in connection with the analysis of ultra high energy cosmic rays.
\end{abstract}

\bigskip


\medskip


\section{Introduction}

We consider the problem of estimating the density $f$ of an independent sample
of points
$X_1,\ldots,X_n$ observed on the $d$-dimensional sphere $\bS^d$ of $\bR^{d+1}$.
Obviously, the most immediate examples of applications appear in
the case $d=2$. However, no major differences arise from considering
the general case.

There is an abundant literature about this type of problems. In
particular, minimax $L^2$ results have been obtained
(see \cite{klemelamini},\cite{klemelalower}). These procedures are generally
obtained using either kernel methods (but in this case the manifold
structure of the sphere is not well taken into account), or using
orthogonal series methods associated with spherical harmonics (and
in this case the 'local performances of the estimator are quite
poor, since spherical harmonics are spread all over the sphere).

In our approach we  focus on two
important points. We aim at a procedure of estimation which
is efficient from a $L^2$ point of view
(as it is a tradition in statistics to evaluate the procedure
with the mean square error). On the other hand, we would like it
to perform satisfactorily also  from a local point of view
(in infinity norm, for instance). To have these two requirements
together seems to us a warrant to have good results in practice. In
effect, it is very difficult to produce a loss function which
reflects at the same time the requirement of clearly seeing the bumps
of the density, of being able to well estimate different
level sets, of testing whether there is a difference between
the northern and southern hemispheres and so on.

In addition, we require this
procedure to be simple to implement, as well as adaptive to
inhomogeneous smoothness.

This type of requirements is generally
well handled using thresholding estimates associated to wavelets.
The  problem requires a special construction adapted to the sphere,
since usual tensorized wavelets will never reflect the manifold
structure of the sphere and will necessarily create unwanted
artifacts. Recently in (\cite{pnarco},\cite{NPW}) a tight frame (i.e. a redundant family) was
produced
which enjoys enough properties to be
successfully  used for density estimation.

The fundamental properties of wavelets are  their concentration in
the Fourier domain as well as in the space domain. Here,
obviously the 'space' domain is the sphere itself whereas the
Fourier domain is now obtained by replacing the 'Fourier' basis by
the basis of Spherical Harmonics which plays an analogous role on
the sphere.

The construction \cite{pnarco},\cite{NPW} produces a family of functions which very
much resemble to wavelets, the needlets, and in particular have very
good concentration properties.

We use these needlets to construct an estimation procedure, and prove
that this procedure attains optimal rates over various spaces of
regularity.

Again, the problem of choosing appropriated spaces of regularity on
the sphere in a serious question, and we decided to consider the
spaces which may be the closest to our natural intuition: those
which generalize to the sphere case the classical H\"older spaces.

In the first section we present (\cite{pnarco}) needlets, and describe spaces of regularity on the sphere.
In the second one we define our estimation procedure, and describe its properties.

The novelties of this paper lie in the application of thresholding to
the needlet coefficients, which gives a very simple and adaptive
procedure which works on the sphere. We also focus here on giving
the results in $\bL_\infty$ norm, and obtain the rates of convergence for many
other loss functions  as a consequence of the previous ones.

Our results are motivated by many recent developments in the area
of observational astrophysics. As an example, we refer to experiments measuring
incoming directions of Ultra High Energy Cosmic Rays, such as the AUGER Observatory
(\verb@http://www.auger.org@).
Here, efficient estimation of the density function of these directional data may yield
crucial insights into the physical mechanisms generating the observations.
More precisely, a uniform density would suggest
the High Energy Cosmic Rays are generated by cosmological effects,
such as the decay of massive particles generated during the Big Bang;
on the other hand,  if these Cosmic Rays are generated by astrophysical phenomena
(such as acceleration into Active Galactic Nuclei),
then we should observe a density function which is highly non-uniform
and tightly correlated with the local distribution of nearby Galaxies.
Massive amount of data in this area are expected to be available in the next few years.
The Auger observatory will be based on two arrays of detectors;
the first one covers an area larger than 3000 Km$^2$ in
Pampa Amarilla (Argentina), and has already started to collect observations:
some preliminary evidence was provided in \cite{collaboration-2008-29}, and
a non-uniform distribution seems to be favored. The whole celestial sphere will actually be covered only when the construction of the northern hemisphere array, due to be built in eastern Colorado,
will be completed, a few years from now. Hence, in the immediate future
efficient statistical techniques will be eagerly requested
for the analysis of the forthcoming datasets.

A survey of statistical methodologies dealing with directional data on the sphere may be found in
\cite{Mardia}, \cite{Jupp}, \cite {Mardia-jupp}.
The generalization of estimation using orthogonal series methods to the case of compact Riemannian manifold can be found in
\cite{hendriks}. See related works in \cite{Kim}, \cite{Ruymgaart}, \cite{hendriks-ruymgaart}, \cite{Pelletier}, \cite{jupp08}.
Kernel methods on the sphere have been investigated in \cite{hallsphere}.
Minimax rates for the equivalent of Sobolev spaces on the sphere associated can be found in \cite{klemelamini}, \cite{klemela},
\cite{klemelalower}.

The plan of the paper is as follows. In \S2 and \S3 we review some background
material on needlets and Besov spaces. \S4 introduces our thresholding estimator,
whose minimax performances are stated in \S5. \S6 shows the performance of the estimators
on some simulated data. \S7--\S9 contain the proofs.
\section{Needlets}\label{ssec-littlewood}
This construction is due to Narcowich, Petrushev and Ward
\cite{pnarco}. Its aim is essentially to build a very well
localized tight frame constructed using spherical harmonics, as
discussed below. It was recently extended with fruitful statistical applications
to more general Euclidean settings (see
\cite{needvd}) and already exploited for estimation and testing problems
in \cite{sphere}, \cite{voronoi}.

Let us denote by $\mathbb{S}^d$, the unit sphere of $\mathbb{R}^{d+1}$.
We denote $dx$ the surface measure of $\mathbb{S}^d$, that is the
unique positive measure on $\bS^d$ which is
invariant by rotation and has total mass $\omega_d=2\pi^{(d+1)/2}/\Gamma(\frac {d+1}2)$.
The following decomposition is well known.
\begin{equation}
\bL^2(dx)=\bigoplus_{l=0}^\infty\cH_l\ ,\label{decomp}
\end{equation}
where $\cH_l$ is the restriction to $\bS^d$ of the homogeneous
polynomials on $\bR^{d+1}$ of degree $l$ which are harmonic (i.e.
$\Delta P=0$, where $\Delta $ is the Laplacian on $\bR^{d+1}$). This
space is called the space of spherical harmonics of degree $l$ (see
\cite{STW}, chap.4, \cite{MR1022665} chap.5). Its dimension is
equal to $g_{l,d}={l+d\choose d}-{l+d-2\choose d}$ and is therefore
of order $l^{d-1}$. The orthogonal projector on $\cH_l$ is given by the
kernel operator
\begin{equation}
\forall \; f\;\in \bL^2(dx),\; P_{\cH_l}f(x)=\int_{\bS^d}L_l(\langle
x,y\rangle) f(y)\,dy
\end{equation}
where $\langle x,y\rangle$ is the
 standard scalar product of $\bR^{d+1}$, and $L_l$ is the Gegenbauer
polynomial with parameter $\frac {d-1}2$ of degree $l$, defined on $[-1,+1]$ and
normalized so that
\begin{equation}\label{lnorm}
\begin{array}{c}
\displaystyle\int_{-1}^1L_l(t)L_k(t) (1-t^2)^{\frac d2-1}\,dt=\frac{g_{l,d}2^{d/2}\Gamma(\frac d2)^2}{\Gamma(d)\,\omega_d^2}\delta_{l,k}=\\
\displaystyle= \frac{g_{l,d}2^{d/2}\Gamma(\frac{d+1}2)^2\Gamma(\frac d2)^2}{4\Gamma(d)\pi^{d+1}}\delta_{l,k}\ .
\end{array}
\end{equation}
For the main situation of interest, $d=2$, the right hand side above is equal to
$\frac {2l+1}{8\pi^2}$. Recall that if $d=2$, the usual normalization of the Legendre polynomial ($L_l(1)=1$) gives the square of their $L^2$ norm equal to
$\frac 2{2l+1}$. Therefore these must be multiplied by $(2l+1)/(4\pi)$, in order to satisfy (\ref{lnorm}).

Let us point
out the following reproducing property of the projection operators:
\begin{equation}
\int_{\bS^d} L_{l}(\langle x,y\rangle)L_{k}(\langle
y,z\rangle)\,dy=\delta_{l,k}L_{l}(\langle x,z\rangle)\ .  \label{auto}
\end{equation}
The construction of needlets is based on the classical Littlewood-Paley
decomposition and a subsequent discretization.

%
Let $\phi $ be a $C^{\infty }$ function on $\mathbb R$, symmetric
and decreasing on $\mathbb R^+$ supported in $|\xi |\leq 1,$ such
that $1\geq \phi (\xi )\geq 0$ and $\phi (\xi )=1$ if $|\xi |\leq
\frac{1}{2} $. We set
\begin{equation*}
b^{2}(\xi )=\phi (\tfrac\xi 2)-\phi (\xi )\geq 0
\end{equation*}%
so that
\begin{equation}
\forall |\xi |\geq 1,\;\qquad \sum_{j\ge 0}b^{2}(\tfrac{\xi }{2^{j}})=1%
\text{ .}  \label{1}
\end{equation}
%
Remark that $b(\xi)\not=0$ only if $\tfrac 12\le|\xi|\le 2$.
Let us now define the operator
$
\Lambda _{j}=\sum_{l\geq 0}b^{2}(\tfrac{l}{2^{j}})L_{l}
$
and the associated kernel
\begin{equation*}
\Lambda _{j}(x,y)=\sum_{l\geq 0}b^{2}(\tfrac{l}{2^{j}})L_{l}(\langle
x,y\rangle)=%
\sum_{2^{j-1}<l<2^{j+1}}b^{2}(\tfrac{l}{2^{j}})L_{l}(\langle
x,y\rangle)\ .
\end{equation*}%
The following proposition is obvious:
\begin{proposition} For every $f\in L^2$
\begin{equation}
f=\lim_{J\rightarrow \infty
}L_{0}(f)+\sum_{j=0}^{J}\Lambda _{j}(f)\ .  \label{rep}
\end{equation}%
Moreover, if $M_{j}(x,y)=\sum_{l\geq
0}b(\frac{l}{2^{j}})L_{l}(\langle x,y\rangle)$, then
\begin{equation}
\Lambda
_{j}(x,y)=\int M_{j}(x,z)M_{j}(z,y)\,dz\ .  \label{sqrt}
\end{equation}
\end{proposition}
%
%
Let
\begin{equation*}
\mathscr{P}_{l}=\bigoplus_{m=0}^{l}\cH_{m}\ ,
\end{equation*}
the space of the restrictions to $\bS^d$ of the polynomials of
degree $\le l$. The following  quadrature formula is true: for all $l\in\mathbb{N}$
there exists a finite subset $\mathscr{X}_{l}\subset\bS^d$ and positive real numbers $%
\lambda_{\eta}>0$, indexed by the elements $\eta\in\mathscr{X}_{l},$
such that
\begin{equation}
\forall f\in\mathscr{P}_{l},\qquad\int_{\bS^d} f(x)\,dx=\sum_{\eta\in\mathscr{X}%
_{l}}\lambda_{\eta}f(\eta)\ .  \label{quadr}
\end{equation}
Then the operator $M_{j}$ defined in the subsection above is such that:
$z\mapsto M_{j}(x,z)\in\mathscr{P}_{[2^{j+1}]}$,
so that
\begin{equation*}
z\mapsto M_{j}(x,z)M_{j}(z,y)\in\mathscr{P}_{[2^{j+2}]}\ ,
\end{equation*}
and we can write:
\begin{equation*}
\Lambda_{j}(x,y)=\int M_{j}(x,z)M_{j}(z,y)dz=\sum_{\eta\in\mathscr{X}%
_{[2^{j+2}]}}\lambda_{\eta}M_{j}(x,\eta)M_{j}(\eta,y)\ .
\end{equation*}
This implies:
\begin{align*}
\Lambda_{j}f(x) & =\int\Lambda_{j}(x,y)f(y)dy=\int\sum_{\eta \in%
\mathscr{X}_{[2^{j+2}]}}\lambda_{\eta}M_{j}(x,\eta)M_{j}(\eta ,y)f(y)dy
\\
& =\sum_{\eta\in\mathscr{X}_{[2^{j+2}]}}\sqrt{\lambda_{\eta}}%
M_{j}(x,\eta)\int {\sqrt{\lambda_{\eta}}M_{j}(y,\eta)}f(y)dy\ .
\end{align*}
We denote
\begin{equation*}
\mathscr{X}_{[2^{j+2}]}=\mathscr{Z}_j,\qquad\psi_{j,\eta}(x):=
\sqrt{\lambda_{\eta}}\,M_{j}(x,\eta)\mbox{ for }\eta\in\mathscr {Z}_j\ .
\end{equation*}
The choice of the sets $\ccl Z_j$ of cubature points is not unique, but one
can impose the conditions
\begin{equation}\label{card}
\tfrac 1c \, 2^{dj}\le \# \mathscr {Z}_j\le c\, 2^{dj},
\qquad \tfrac 1c 2^{-dj}\le\lambda_\eta\le c 2^{-dj}
\end{equation}
for some $c>0$. Actually in the simulations of \S\ref{simulations} we make use
of some sets of cubature points for $d=2$ such that $\# \mathscr {Z}_j=2^{2j+4}$
exactly (the corresponding weights being however not identical).
We have,  using (\ref{rep})
\begin{equation}\label{develop}
f =  L_{0}(f) + \sum_{j}\sum_{\eta\in \mathscr{Z}_j}
\langle f,\psi_{j,\eta}\rangle_{\bL^2(\bS^d)}\psi_{j,\eta}\ .
\end{equation}
\begin{figure}[h]
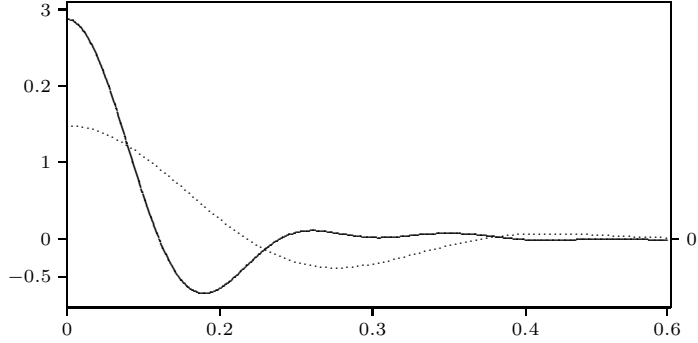

\hbox to\hsize\bgroup\hss
\beginpicture
\setcoordinatesystem units <4truein,.4truein> \setplotarea x from 0
to 0.79, y from -0.9 to 3.1 \axis bottom  ticks short withvalues
{$\scriptstyle0$}  {$\scriptstyle0.2$} {$\scriptstyle0.3$}
{$\scriptstyle0.4$}  {$\scriptstyle0.6$} {$\scriptstyle\frac \pi4$}
/ at 0  .2  .4  .6  0.7853 / / \axis left shiftedto x=0 ticks short
withvalues {$\scriptstyle-0.5$} {$\scriptstyle0$} {$\scriptstyle1$}
{$\scriptstyle0.2$} {$\scriptstyle3$} / at -0.5 0 1 2 3 / / \axis
right ticks short withvalues $\scriptstyle0$ / at 0 / / \axis top /
\setdots <2pt> \plot
   0.           1.4759783
    0.0078540    1.4733059
    0.0157080    1.4653061
    0.0235619    1.4520311
    0.0314159    1.4335673
    0.0392699    1.4100347
    0.0471239    1.3815858
    0.0549779    1.3484047
    0.0628319    1.310705
    0.0706858    1.2687287
    0.0785398    1.2227437
    0.0863938    1.1730421
    0.0942478    1.1199375
    0.1021018    1.0637624
    0.1099557    1.0048661
    0.1178097    0.9436115
    0.1256637    0.8803720
    0.1335177    0.8155293
    0.1413717    0.7494698
    0.1492257    0.6825818
    0.1570796    0.6152528
    0.1649336    0.5478663
    0.1727876    0.4807989
    0.1806416    0.4144179
    0.1884956    0.3490782
    0.1963495    0.2851202
    0.2042035    0.2228674
    0.2120575    0.1626239
    0.2199115    0.1046728
    0.2277655    0.0492743
    0.2356194  -0.0033359
    0.2434734  -0.0529478
    0.2513274  -0.0993780
    0.2591814  -0.1424710
    0.2670354  -0.1820992
    0.2748894  -0.2181636
    0.2827433  -0.2505937
    0.2905973  -0.2793472
    0.2984513  -0.3044100
    0.3063053  -0.3257948
    0.3141593  -0.3435411
    0.3220132  -0.3577132
    0.3298672  -0.3683999
    0.3377212  -0.3757119
    0.3455752  -0.3797814
    0.3534292  -0.3807592
    0.3612832  -0.3788135
    0.3691371  -0.3741279
    0.3769911  -0.3668988
    0.3848451  -0.3573340
    0.3926991  -0.3456499
    0.4005531  -0.3320696
    0.4084070  -0.3168208
    0.4162610  -0.3001334
    0.4241150  -0.2822376
    0.4319690  -0.2633619
    0.4398230  -0.2437308
    0.4476770  -0.2235635
    0.4555309  -0.2030715
    0.4633849  -0.1824576
    0.4712389  -0.1619141
    0.4790929  -0.1416217
    0.4869469  -0.1217479
    0.4948008  -0.1024466
    0.5026548  -0.0838570
    0.5105088  -0.0661029
    0.5183628  -0.0492924
    0.5262168  -0.0335176
    0.5340708  -0.0188544
    0.5419247  -0.0053624
    0.5497787    0.0069143
    0.5576327    0.0179473
    0.5654867    0.0277233
    0.5733407    0.0362432
    0.5811946    0.0435213
    0.5890486    0.0495849
    0.5969026    0.0544728
    0.6047566    0.0582344
    0.6126106    0.0609286
    0.6204645    0.0626226
    0.6283185    0.0633907
    0.6361725    0.0633131
    0.6440265    0.0624743
    0.6518805    0.0609624
    0.6597345    0.0588675
    0.6675884    0.0562804
    0.6754424    0.0532918
    0.6832964    0.0499912
    0.6911504    0.0464653
    0.6990044    0.0427978
    0.7068583    0.0390681
    0.7147123    0.0353505
    0.7225663    0.0317137
    0.7304203    0.0282202
    0.7382743    0.0249257
    0.7461283    0.0218790
    0.7539822    0.0191214
    0.7618362    0.0166867
    0.7696902    0.0146012
    0.7775442    0.0128837
    0.7853982    0.0115456
/
\setsolid
\plot
          0.           2.8755571
    0.0078540    2.8556175
    0.0157080    2.7963007
    0.0235619    2.6990972
    0.0314159    2.5664381
    0.0392699    2.4016188
    0.0471239    2.2086943
    0.0549779    1.9923529
    0.0628319    1.75777
    0.0706858    1.5104499
    0.0785398    1.256059
    0.0863938    1.0002576
    0.0942478    0.7485351
    0.1021018    0.5060558
    0.1099557    0.2775185
    0.1178097    0.0670357
    0.1256637  -0.1219645
    0.1335177  -0.2868107
    0.1413717  -0.4256329
    0.1492257  -0.5373800
    0.1570796  -0.6218081
    0.1649336  -0.6794442
    0.1727876  -0.7115248
    0.1806416  -0.7199139
    0.1884956  -0.7070040
    0.1963495  -0.6756049
    0.2042035  -0.6288240
    0.2120575  -0.569943
    0.2199115  -0.5022971
    0.2277655  -0.4291590
    0.2356194  -0.3536330
    0.2434734  -0.2785628
    0.2513274  -0.2064545
    0.2591814  -0.1394186
    0.2670354  -0.0791300
    0.2748894  -0.0268082
    0.2827433    0.0167833
    0.2905973    0.0513209
    0.2984513    0.0768879
    0.3063053    0.0939279
    0.3141593    0.1031888
    0.3220132    0.1056576
    0.3298672    0.1024922
    0.3377212    0.0949514
    0.3455752    0.0843274
    0.3534292    0.0718823
    0.3612832    0.0587924
    0.3691371    0.0461016
    0.3769911    0.0346846
    0.3848451    0.0252222
    0.3926991    0.0181875
    0.4005531    0.0138435
    0.4084070    0.0122517
    0.4162610    0.0132891
    0.4241150    0.0166737
    0.4319690    0.0219951
    0.4398230    0.0287498
    0.4476770    0.0363778
    0.4555309    0.0442990
    0.4633849    0.0519485
    0.4712389    0.0588073
    0.4790929    0.0644288
    0.4869469    0.0684590
    0.4948008    0.0706506
    0.5026548    0.0708702
    0.5105088    0.0690983
    0.5183628    0.0654241
    0.5262168    0.0600338
    0.5340708    0.0531950
    0.5419247    0.0452375
    0.5497787    0.0365318
    0.5576327    0.0274668
    0.5654867    0.0184282
    0.5733407    0.0097780
    0.5811946    0.0018368
    0.5890486  -0.0051308
    0.5969026  -0.0109268
    0.6047566  -0.0154259
    0.6126106  -0.0185777
    0.6204645  -0.0204044
    0.6283185  -0.0209946
    0.6361725  -0.0204946
    0.6440265  -0.0190956
    0.6518805  -0.0170209
    0.6597345  -0.0145108
    0.6675884  -0.0118084
    0.6754424  -0.0091452
    0.6832964  -0.0067293
    0.6911504  -0.0047344
    0.6990044  -0.0032925
    0.7068583  -0.0024887
    0.7147123  -0.0023594
    0.7225663  -0.0028927
    0.7304203  -0.0040327
    0.7382743  -0.0056847
    0.7461283  -0.0077237
    0.7539822  -0.0100031
    0.7618362  -0.0123649
    0.7696902  -0.0146500
    0.7775442  -0.0167072
    0.7853982  -0.0184027
/

\endpicture
\hss\egroup \caption{\footnotesize  The value of the needlet $\psi_{j,\xi}$  as a function of the distance
from the cubature point $\xi$ for $j=3$ (dots) and $j=4$ (solid).
This emphasizes the localization properties
of needlets as $j$ increases.\label{needlet}}
\end{figure}
The main result of Narcowich, Petrushev and
Ward, \cite{pnarco}  is the following localization property of the
 $\psi_{j,\eta}$, called needlets: for any $k$ there exists
 a constant $c_{k}$ such that, for every $\xi\in
\mathbb{S}^d$:
\begin{equation}\label{localized}
|\psi _{j,\eta }(\xi )|\leq \frac{c_{k}2^{jd/2}}{(1+2^{jd/2}d(\eta
,\xi))^{k}}\ \cdotp
\end{equation}
where $d$ is the natural geodesic distance on the sphere  (for $d=2$,  $d(\xi,\eta)= \arccos \langle\eta ,\xi \rangle$
). In other words needlets
are almost exponentially localized around any cubature point, which
motivates their name.
%
From this localization property it follows (see \cite{pnarco}) that for $1\le p\le+\infty$
there exists positive constants $c_p,C_p$ such that
\begin{equation}\label{lp}
c_p2^{jd({\frac 12}-\frac 1p)}\le \|\psijk\|_p\le C_p2^{jd({\frac 12}-\frac 1p)}\ .
\end{equation}
Also the following holds
\begin{lemma}\label{lemma-inequality} 1) For every $0< p\le+\infty$
\begin{equation}\label{inequality}
\Big\Vert \sum_{\xi\in \mathscr{Z}_j}\lambda_\xi \psi_{j,\xi}\Big\Vert_p\le
c 2^{jd(\frac12-\frac 1p)}\Bigl(\sum_{\xi\in \mathscr{Z}_j}
|\lambda_\xi|^p\Bigr)^{1/p}\ .
\end{equation}
2) For every $1\le p\le+\infty$ \begin{equation}\label{inequality2}
\Big( \sum_{\xi\in \mathscr{Z}_j}|\langle f,\psi_{j\xi}\rangle|^p\Big)^{1/p}
2^{jd(\frac12-\frac 1p)}\le
c\Vert f\Vert_p\ .
\end{equation}
\end{lemma}
\proof Let us prove \tref{inequality} for $p=+\infty$. Using \tref{localized} and
Lemma 6 of \cite{sphere}
$$
\displaylines{
\sup_{x\in\mathbb{S}^d} \Big| \sum_{\xi\in \mathscr{Z}_j}\lambda_\xi \psi_{j,\xi}\Big|
\le
\sup_{\xi\in\ccl Z_j}
|\lambda_\xi |\sup_{x\in\mathbb{S}^d} \sum_{\xi\in \mathscr{Z}_j}|\psi_{j,\xi}(x)|\le\cr
\le
\sup_{\xi\in\ccl Z_j} |\lambda_\xi |c_3\sup_{x\in\mathbb{S}^d}
\sum_{\xi\in \mathscr{Z}_j}\frac{2^{jd/2}}{(1+2^{jd/2}d(\xi,x))^3}\le c_12^{jd/2}
\sup_{\xi\in\ccl Z_j} |\lambda_\xi |\ .\cr
}
$$
If $1\le p<+\infty$, by H\"older inequality, if $\frac 1p+\frac 1{p'}=1$ so that
$\frac p{p'}=p-1$,
$$
\displaylines{
\Bigl(\sum_{\xi\in \mathscr{Z}_j}|\lambda_\xi \psi_{j,\xi}(x)|\Bigr)^p=
\Bigl(\sum_{\xi\in \mathscr{Z}_j}|\lambda_\xi| |\psi_{j,\xi}(x)|^{1/p}
|\psi_{j,\xi}(x)|^{1/p'}\Big)^p\le\cr
\le\Big( \sum_{\xi\in \mathscr{Z}_j}|\lambda_\xi| ^p|\psi_{j,\xi}(x)|
\Big)\Big( \sum_{\xi\in \mathscr{Z}_j} |\psi_{j,\xi}(x)|
\Big)^{p-1}\le\cr
\le c_32^{j\frac d2\,(p-1)}\sum_{\xi\in \mathscr{Z}_j}|\lambda_\xi| ^p|\psi_{j,\xi}(x)|
}
$$
where the last inequality comes again from \tref{localized} and
Lemma 6 of \cite{sphere}.
Now integrating and using \tref{lp} for $p=1$,
$$
\Big\Vert \sum_{\xi\in \mathscr{Z}_j}\lambda_\xi \psi_{j,\xi}(x)\Big\Vert_p^p\le
2^{j\frac d2\,(p-1)}\sum_{\xi\in \mathscr{Z}_j}|\lambda_\xi|^p \|\psi_{j,\xi}\|_1
\le 2^{j\frac d2(p-2)}\sum_{\xi\in \mathscr{Z}_j}|\lambda_\xi|^p
$$
from which \tref{inequality} follows. The remaining case $0< p\le1$ follows immediately by
subadditivity, as
$$
\Big\Vert \sum_{\xi\in \mathscr{Z}_j}\lambda_\xi \psi_{j,\xi}(x)\Big\Vert_p^p
\le \sum_{\xi\in \mathscr{Z}_j}|\lambda_\xi|^p \|\psi_{j,\xi}(x)\|_p^p\ .
$$
As for 2) clearly if $p=+\infty$
$$
\displaylines{
C 2^{jd/2}\sup_{\xi  \in \ccl Z_j}|  \langle f, \phi_{j,\xi} \rangle |
\leq C 2^{jd/2}\sup_{\xi  \in \ccl Z_j}  \int |  f (x) | | \phi_{j,\xi} (x) |\, dx   \le\cr
\leq C 2^{jd/2}\| f \|_\infty
 \sup_{\xi  \in \ccl Z_j } \|\phi_{j,\xi} \|_1   \leq C'  \| f \|_\infty\cr
 }
$$
and if $p=1$
$$\displaylines{
\sum_{\xi  \in \ccl Z_j }|  \langle f, \phi_{j,\xi} \rangle | 2^{-jd/2}
\leq 2^{-jd/2}\sum_{\xi  \in \ccl Z_j}  \int |  f(x)| | \phi_{j,\xi} (x) |\, dx   =\cr
= 2^{-jd/2}\int |  f(x)| \sum_{\xi  \in \ccl Z_j} | \phi_{j,\xi} (x) |  \,dx
\leq C \|f\|_1\ .\cr
}
$$
Let now $1<p<\infty $    
$$
\displaylines{
\sum_{\xi  \in \ccl Z_j}|  \langle f, \phi_{j,\xi} \rangle |^p 2^{jd(p/2-1)} \leq
2^{jd(p/2-1)}\sum_{\xi  \in \ccl Z_j} \Big( \int |  f (x)|  | \phi_{j,\xi} (x) |
\, dx\Big)^p \ .\cr}
$$
But, by Holder inequality, for $p'$ such that $\frac 1p + \frac 1p' =1, then \frac p{p'} =p-1$ and
$$
\displaylines{
\Big( \int |  f (x)|  | \phi_{j,\xi} (x) | dx\Big)^p= \Big( \int |  f (x)|
| \phi_{j,\xi} (x) |^{1/p} | \phi_{j,\xi} (x) |^{1/p'}\, dx\Big)^p\le\cr
\leq  \int |  f (x)|^p  | \phi_{j,\xi} (x) |\,  dx
\Big(\int | \phi_{j,\xi} (x) |\, dx \Big)^{p-1} \!\! \!\!=
 \int |  f (x)|^p  | \phi_{j,\xi} (x) |\,  dx \| \phi_{j,\xi} \|_1^{p-1} .\cr
 }
$$
So
$$
\displaylines{
\sum_{\xi  \in \ccl Z_j}|  \langle f, \phi_{j,\xi} \rangle |^p \| 2^{jd(p/2-1)}
\leq  2^{-j\frac d2(p-1)}    2^{jd(p/2-1)}\! \!\sum_{\xi  \in \ccl Z_j}
\int |  f (x)  |^p| \phi_{j,\xi} (x) |\, dx=\cr
=
2^{-jd/2}\int |  f (x) |^p     \sum_{\xi  \in \ccl Z_j}
|\phi_{j,\xi} (x) |dx  \leq C\|f\|_p^p\ .\cr
}
$$
%

\begin{example}\label{l2norm-of-needlet}\rm
Relation (\ref{lp}) for $p=2$ states that the $L^2$ norm of
$\psi_{j\xi}$ is bounded with respect to $j$ and also bounded away
from $0$ from below. Assume $d=2$. Then using (\ref{auto}) it is actually easy to see
that, keeping in mind that $L_l(1)=\frac {2l+1}{4\pi}$,
$$
\Vert \psi_{j\xi}\Vert^2_2=\lambda_\xi\sum_{l\geq
0}b^2(\tfrac{l}{2^{j}})L_{l}(1)=\frac {\lambda_\eta}{4\pi}\sum_{l\geq
0}b^2(\tfrac{l}{2^{j}})(2l+1)\ .
$$
Assuming that the cubature points are of cardinality $2^{2j+4}$
and that they sum up to $4\pi$, $\lambda_\eta\sim 4\pi\cdot
2^{-2j-4}$ as $j\to\infty$. If the previous relation were an equality we could recognize
in the right hand term the Riemann sum
$$
\frac 18\Bigl(\tfrac{1}{2^{j}}\sum_{l\geq
0}b^2(\tfrac{l}{2^{j}})\tfrac{l}{2^{j}}+\tfrac{1}{2^{2j}}\sum_{l\geq
0}b^2(\tfrac{l}{2^{j}})\Bigr)
$$
that converges, as $j\to\infty$, to the integral
$$
I=\frac18\int_{1/2}^2 tb^2(t)\, dt
$$
which depends on the choice of the function $b$. This $L^2$ norm shall appear
frequently in the sequel. For instance, if we write down the development \tref{develop} or
the function $f=\psi_{j_0\xi_0}$, then the coefficient
$\beta_{j_0\xi_0}=\langle f,\psi_{j,\eta}\rangle$ would be exactly equal to
$\Vert \psi_{j\xi}\Vert^2_2$. As it is clear that it would be desirable for
this coefficient to be as large as possible, the value of the integral above can be seen
as a measure of the localization properties of the system of needlets and can be used as
a criterion of goodness of the choice of the function $b$. With the choice we made
(see \S\ref{simulations}) the quantity $I$ above is $\simeq 0.107$.
\end{example}
\section{Besov spaces on the sphere and needlets}
In this section we summarize the main properties of Besov spaces and needlets, as
established in \cite{pnarco}.

Let $f:\mathbb{S}^d\to\R$ a measurable function. We define
$$
E_k(f,r)=\inf_{P\in\mathscr{P}_k}\Vert f-P\Vert_r
$$
the infimum of the distances in $L^r$ of $f$ from the polynomials of degree $k$.
Then the Besov space $B^s_{r, q}$ is defined as the space of functions such that
$$
f\in L^r\mbox{ and } \Bigl(\sum_{k=0}^\infty (k^s E_k(f,r))^q\tfrac 1k\Bigr)^{1/q}<+\infty\ .
$$
Remarking that $k\to E_k(f,r))$ is decreasing, by a standard condensation argument this is
equivalent to
$$
f\in L^r\mbox{ and } \Bigl(\sum_{j=0}^\infty (2^{js} E_{2^j}(f,r))^q\Bigr)^{1/q}<+\infty\ .
$$
\begin{theorem}
Let $1\le r\le+\infty$, $s>0$, $0\le q\le +\infty$. Let $f$ a measurable function
and define
$$
\langle f,\psi_{j,\xi}\rangle=\int_{\mathbb{S}^d}f(x) \psi_{j,\xi}(x)\, dx\enspace\mathop{=}^{def}
\enspace \beta_{j,\xi}
$$
provided the integrals exists. Then $f\in B^s_{r, q}$ if and only if, for every
$j=1,2,\dots$,
$$
\Bigl(\sum_{\xi\in\mathscr{X}_j}
(\beta_{j,\xi}\Vert\psi_{j,\xi}\Vert_r)^r\Bigr)^{1/r}
=2^{-js}\delta_j
$$
where $(\delta_j)_j\in\ell_q$.
\end{theorem}
As
$$
c2^{jd(\frac12-\frac 1r)}\le
\Vert\psi_{j,\xi}\Vert_r\le C 2^{jd(\frac12-\frac 1r)}
$$
for some positive constants $c,C$, the Besov space $B^s_{r, q}$ turns out to be
a Banach space associate to the norm
\begin{equation}\label{bnorm}
\displaystyle\|f\|_{B^s_{r, q}}
:= \|(2^{j[s+d(\frac 12-\frac 1r)]}\|(\beta_{j\eta})_{\eta\in\Z_j}
\|_{\ell_r})_{j\ge 0}\|_{\ell_q}<\infty\ .\\
\end{equation}
In the sequel we shall denote by $B^s_{r,q}(M)$ the ball of radius $M$ of the
Besov space $B^s_{r,q}$.

\begin{theorem}\label{Besov-embedding} (The Besov embedding)
If $ p\leq r \leq \infty$ then $B^{s}_{r,q} \subseteq B^{s}_{p,q}$.
If $s >d(\frac 1r - \frac 1p)$,
$$
r  \leq p \leq \infty  \Rightarrow  B^{s}_{r,q} \subseteq B^{s-d(1/r - 1/p)}_{p,q}\ .
$$
\end{theorem}
\noindent{\bf Proof} By hypothesis
$$
2^{jd(\frac 12-\frac 1r)}   \Big( \sum_{\xi \in \ccl Z_j}  |\beta_{j,\xi}|^r
\Big)^{1/r}  \leq  \delta_j 2^{-js},\qquad\hbox to0pt{$(\delta_j)_j\in \ell_q$}\ .
$$
Let $p\leq r \leq \infty$, then
$$
\displaylines{
2^{jd( \frac 12- \frac 1p)}\Big( \sum_{\xi \in \ccl Z_j}  |\beta_{j,\xi}|^p \Big)^{1/p}=\cr
= 2^{jd/2}  (2^{-jd} {\rm card} (\ccl Z_j))^{1/p}  \Big( \frac 1{{\rm card} (\ccl Z_j)}
\sum_{\xi \in \ccl Z_j}  |\beta_{j,\xi}|^p \Big)^{1/p}\le\cr
\leq  2^{jd/2}  (2^{-jd} {\rm card} (\ccl Z_j))^{1/p}
\Big( \frac 1{{\rm card}(\ccl Z_j)} \sum_{\xi \in \ccl Z_j}  |\beta_{j,\xi}|^r \Big)^{1/r}=\cr
=  2^{jd(\frac 12-\frac 1r)}  (2^{-jd} {\rm card}(\ccl Z_j))^{\frac 1p-\frac 1r}
\Big( \sum_{\xi \in \ccl Z_j}  |\beta_{j,\xi}|^r \Big)^{1/r}\le\cr
\leq C
2^{jd(\frac 12-\frac 1r)}   \Big( \sum_{\xi \in \ccl Z_j}  |\beta_{j,\xi}|^r \Big)^{1/r}
\leq C  \delta_j 2^{-js}\ .\cr
}
$$
On the other hand, if $ r  \leq p \leq \infty  $,
$$
\displaylines{
2^{jd (\frac 12- \frac 1p)}    \Big( \sum_{\xi \in \ccl Z_j}  |\beta_{j,\xi}|^p \Big)^{1/p}
\leq   2^{jd (\frac 12- \frac 1p)}    \Big( \sum_{\xi \in \ccl Z_j}  |\beta_{j,\xi}|^r \Big)^{1/r}=\cr
= 2^{jd (\frac 1r- \frac 1p)} 2^{jd (\frac 12- \frac 1r)}
\Big( \sum_{\xi \in \ccl Z_j}  |\beta_{j,\xi}|^r \Big)^{1/r}\le\cr
\leq 2^{jd (\frac 1r- \frac 1p)}  \delta_j 2^{-js} =  \delta_j 2^{-j(s -d (\frac 1r- \frac 1p))}\ .\cr
}
$$
\section{Needlet estimation of a density on the sphe\-re}
Let us suppose that we observe $X_1,\ldots,X_n$, i.i.d. random
variables taking values on the sphere having common density $f$ with
respect to $dx$.
$f$ can be decomposed using the frame of needlets described above.
$$
f =\frac 1{|\mathbb S^d|}+\sum_{j\ge 0}\sum_\eta
\bjk\psijk\ .
$$
The needlet estimator
is   based on  hard thresholding of a needlet expansion as follows.
We start by letting:
\begin{align}
\label{eq:hatbeta}\hbjk&:=\frac 1n\sum_{i=1}^n\psijk(X_i)\\
\label{eq:need}\hf&:=\frac 1{|\mathbb S^d|}+\sum_{j=0}^J\sum_{\eta\in
\ccl Z_j} \,\hbjk\,\psijk\,{1}_{\{|\hbjk|\geq \,\kappa\,c_n\}}\ .
\end{align}
The  tuning parameters of the needlet estimator are:

$\bullet$ The range $J=J(n)$ of resolution levels (frequencies) where
the approximation (\ref{eq:need}) is used:
$$
\displaylines{
\Lambda_n=\{(j,\eta),\;0\le j\le J,\; \eta\in \Z_j\},\cr
}
$$
We shall see that the choice $2^{J}=\big( \frac n{\log n}\big)^{\frac 1d}$
is appropriate.

$\bullet$ The threshold constant $\kappa$.  
Evaluations of $\kappa$ are
given in the following Section and also discussed in \S\ref{simulations}.

$\bullet$ $c_n$: is a sample size-dependent scaling factor. We shall see that an appropriate
choice is
$$
c_n=\Big( {\log n\over n} \Big)^{1/2}.
$$
\begin{example} \rm In order to give a better intuition about the localization and
near absence of correlation of the needlet coefficients, let us consider the case of a sample
$X_1,\dots,X_n$ of i.i.d. r.v.'s uniform on the sphere $\mathbb{S}^2$ of $\mathbb{R}^3$.
Then the distribution of the r.v. $\langle x,X_i\rangle$ is uniform on
the interval $[-1,1]$ and, if $\widehat \beta_{\eta j}$, $\widehat \beta_{\xi j}$ are the corresponding needlet
coefficients associated to the cubature points $\eta,\xi$ then of course they are
centered r.v.'s and, thanks to (\ref{auto}), their covariance is equal to
$$
E(\widehat \beta_{\eta j}\widehat \beta_{\xi j})=\sqrt{\lambda_\eta\lambda_\xi}\sum_{l\geq
0}b^2(\tfrac{l}{2^{j}})L_{l}(\langle \xi,\eta\rangle)\ .
$$
Setting $\eta=\xi$  we find
$$
\Var(\widehat \beta_{\eta j})=\Vert \psi_{j\xi}\Vert_2^2
$$
which is a quantity already discussed in Example \ref{l2norm-of-needlet}.
As for the cor\-re\-la\-tion between coefficients, it is given by the function
$\th\to\lambda_\eta \sum_{l\geq
0}b^2(\tfrac{l}{2^{j}})L_{l}(\cos\th)$, whose graph, for some values of $j$ is
plotted in Figure \ref{needletcov}.

\begin{figure}[h]
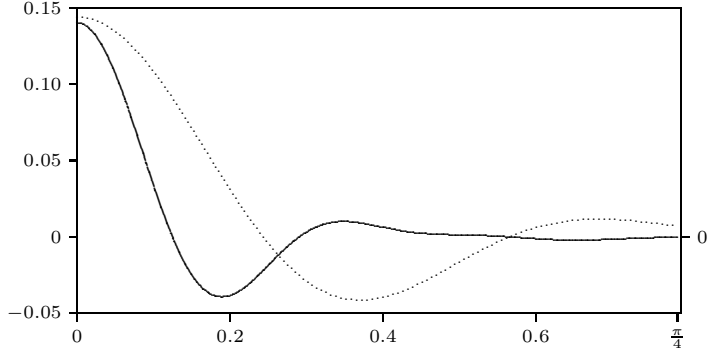

\hbox to\hsize\bgroup\hss
\beginpicture
\setcoordinatesystem units <4truein,8truein>
\setplotarea x from 0 to 0.79, y from -0.05 to 0.15
\axis bottom  ticks short
withvalues {$\scriptstyle0$}  {$\scriptstyle0.2$}  {$\scriptstyle0.4$}
 {$\scriptstyle0.6$}  {$\scriptstyle\frac \pi4$}  / at 0   .2  .4  .6  0.7853 / /
\axis left shiftedto x=0 ticks short withvalues {$\scriptstyle-0.05$} {$\scriptstyle0$}
{$\scriptstyle0.05$} {$\scriptstyle0.10$} {$\scriptstyle0.15$}  / at -0.05 0 0.05 0.10 0.15 / /
\axis right ticks short withvalues $\scriptstyle0$ / at 0 / /
\axis top /
\setdots <2pt>
\plot
    0.           0.1440762
    0.0078540    0.1438386
    0.0157080    0.1431273
    0.0235619    0.1419464
    0.0314159    0.1403026
    0.0392699    0.1382054
    0.0471239    0.1356668
    0.0549779    0.1327014
    0.0628319    0.1293260
    0.0706858    0.1255600
    0.0785398    0.1214245
    0.0863938    0.1169429
    0.0942478    0.1121403
    0.1021018    0.1070434
    0.1099557    0.1016806
    0.1178097    0.0960812
    0.1256637    0.0902758
    0.1335177    0.0842957
    0.1413717    0.0781728
    0.1492257    0.0719396
    0.1570796    0.0656287
    0.1649336    0.0592725
    0.1727876    0.0529033
    0.1806416    0.0465530
    0.1884956    0.0402527
    0.1963495    0.0340327
    0.2042035    0.0279223
    0.2120575    0.0219494
    0.2199115    0.0161409
    0.2277655    0.0105216
    0.2356194    0.0051150
    0.2434734  -0.0000573
    0.2513274  -0.0049756
    0.2591814  -0.0096224
    0.2670354  -0.0139820
    0.2748894  -0.0180412
    0.2827433  -0.0217887
    0.2905973  -0.0252156
    0.2984513  -0.0283152
    0.3063053  -0.0310829
    0.3141593  -0.0335166
    0.3220132  -0.0356160
    0.3298672  -0.0373831
    0.3377212  -0.0388221
    0.3455752  -0.0399387
    0.3534292  -0.0407410
    0.3612832  -0.0412384
    0.3691371  -0.0414421
    0.3769911  -0.0413649
    0.3848451  -0.0410207
    0.3926991  -0.0404250
    0.4005531  -0.0395939
    0.4084070  -0.0385447
    0.4162610  -0.0372955
    0.4241150  -0.0358649
    0.4319690  -0.0342718
    0.4398230  -0.0325356
    0.4476770  -0.0306759
    0.4555309  -0.0287119
    0.4633849  -0.0266630
    0.4712389  -0.0245482
    0.4790929  -0.0223860
    0.4869469  -0.0201944
    0.4948008  -0.0179907
    0.5026548  -0.0157914
    0.5105088  -0.0136122
    0.5183628  -0.0114677
    0.5262168  -0.0093717
    0.5340708  -0.0073367
    0.5419247  -0.0053740
    0.5497787  -0.0034941
    0.5576327  -0.0017058
    0.5654867  -0.0000170
    0.5733407    0.0015657
    0.5811946    0.0030370
    0.5890486    0.0043928
    0.5969026    0.0056302
    0.6047566    0.0067473
    0.6126106    0.0077434
    0.6204645    0.0086189
    0.6283185    0.0093752
    0.6361725    0.0100145
    0.6440265    0.0105399
    0.6518805    0.0109555
    0.6597345    0.0112658
    0.6675884    0.0114761
    0.6754424    0.0115921
    0.6832964    0.0116201
    0.6911504    0.0115666
    0.6990044    0.0114387
    0.7068583    0.0112433
    0.7147123    0.0109878
    0.7225663    0.0106792
    0.7304203    0.0103249
    0.7382743    0.0099320
    0.7461283    0.0095074
    0.7539822    0.0090579
    0.7618362    0.0085900
    0.7696902    0.0081097
    0.7775442    0.0076228
    0.7853982    0.0071347
/
\setsolid
\plot
       0.           0.1403755
    0.0078540    0.1394905
    0.0157080    0.1368551
    0.0235619    0.1325279
    0.0314159    0.1266044
    0.0392699    0.1192150
    0.0471239    0.1105206
    0.0549779    0.1007089
    0.0628319    0.0899887
    0.0706858    0.0785847
    0.0785398    0.0667309
    0.0863938    0.0546652
    0.0942478    0.0426228
    0.1021018    0.0308301
    0.1099557    0.0194997
    0.1178097    0.0088250
    0.1256637  -0.0010240
    0.1335177  -0.0099042
    0.1413717  -0.0177020
    0.1492257  -0.0243348
    0.1570796  -0.0297523
    0.1649336  -0.0339356
    0.1727876  -0.0368965
    0.1806416  -0.0386754
    0.1884956  -0.0393387
    0.1963495  -0.0389754
    0.2042035  -0.0376932
    0.2120575  -0.0356144
    0.2199115  -0.0328717
    0.2277655  -0.0296035
    0.2356194  -0.0259497
    0.2434734  -0.0220478
    0.2513274  -0.0180287
    0.2591814  -0.0140139
    0.2670354  -0.0101124
    0.2748894  -0.0064185
    0.2827433  -0.0030105
    0.2905973    0.0000501
    0.2984513    0.0027196
    0.3063053    0.0049709
    0.3141593    0.0067932
    0.3220132    0.0081906
    0.3298672    0.0091803
    0.3377212    0.0097906
    0.3455752    0.0100587
    0.3534292    0.0100286
    0.3612832    0.0097482
    0.3691371    0.0092679
    0.3769911    0.0086377
    0.3848451    0.0079057
    0.3926991    0.0071166
    0.4005531    0.0063099
    0.4084070    0.0055197
    0.4162610    0.0047732
    0.4241150    0.0040913
    0.4319690    0.0034878
    0.4398230    0.0029702
    0.4476770    0.0025404
    0.4555309    0.0021948
    0.4633849    0.0019258
    0.4712389    0.0017228
    0.4790929    0.0015725
    0.4869469    0.0014609
    0.4948008    0.0013736
    0.5026548    0.0012968
    0.5105088    0.0012182
    0.5183628    0.0011274
    0.5262168    0.0010164
    0.5340708    0.0008798
    0.5419247    0.0007150
    0.5497787    0.0005219
    0.5576327    0.0003029
    0.5654867    0.0000626
    0.5733407  -0.0001928
    0.5811946  -0.0004560
    0.5890486  -0.0007187
    0.5969026  -0.0009727
    0.6047566  -0.0012098
    0.6126106  -0.0014228
    0.6204645  -0.0016052
    0.6283185  -0.0017521
    0.6361725  -0.0018602
    0.6440265  -0.0019277
    0.6518805  -0.0019546
    0.6597345  -0.0019421
    0.6675884  -0.0018934
    0.6754424  -0.0018125
    0.6832964  -0.0017044
    0.6911504  -0.0015748
    0.6990044  -0.0014298
    0.7068583  -0.0012754
    0.7147123  -0.0011175
    0.7225663  -0.0009613
    0.7304203  -0.0008114
    0.7382743  -0.0006715
    0.7461283  -0.0005445
    0.7539822  -0.0004320
    0.7618362  -0.0003349
    0.7696902  -0.0002529
    0.7775442  -0.0001851
    0.7853982  -0.0001299
/
\endpicture
\hss\egroup \caption{\footnotesize The decay of the covariance of $\widehat\beta_{j,\xi}$ and $\widehat\beta_{j,\eta}$ as a function of the distance
between the cubature points $\xi$ and $\eta$ for $j=3$ (dots) and $j=4$ (solid) (case of a uniformly distributed sample).\label{needletcov}}
\end{figure}
\end{example}
\begin{rem}\rm
Whereas coefficients associated to cubature points that are not too close
are only slightly correlated, the random needlet coefficients
$\widehat\beta_{j,\eta}$, $\eta\in\ccl Z_j$ are not independent and they even satisfy
the linear relation
$$
\sum_{\eta\in\ccl Z_j}\sqrt{\lambda_\eta}\,\widehat\beta_{j,\eta}=0\ .
$$
This comes from the fact that, as $y\to L_l(\langle y,x\rangle)$ for $l\le 2^j$ is a
polynomial of degree $\le 2^{2j}$, one has
\begin{equation}\label{zerosum}
 \sum_{\eta\in\ccl Z_j}\lambda_\eta L_l(\langle \eta,x\rangle)=
\int_{\mathbb{S}^d}L_l(\langle y,x\rangle)\, dy=0\ .
\end{equation}
Therefore
$$
\sum_{\eta\in\ccl Z_j}\sqrt{\lambda_\eta}\,\widehat\beta_{j,\eta}=
\frac 1n\sum_{i=1}^n\sum_{l>0}b(\tfrac l{2^j})
\underbrace{\sum_{\eta\in \ccl Z_j} \lambda_\eta L_l(\langle\eta,X_i\rangle)}_{=0}=0\ .
$$
Relation (\ref{zerosum}) also implies that, for a given square integrable function $f$
on $\mathbb{S}^d$, $\sum_{\eta\in\ccl Z_j}\sqrt{\lambda_\eta}\beta_{j,\eta}=0$. Actually
$$
\sum_{\eta\in\ccl Z_j}\sqrt{\lambda_\eta}\beta_{j,\eta}=
\sum_{l>0}b(\tfrac l{2^j})\int_{\mathbb{S}^d}\underbrace{\sum_{\eta\in\ccl Z_j}\lambda_\eta
L_l(\langle \eta,x\rangle)}_{=0}f(x)\, dx\ .
$$
\end{rem}
\section{Minimax rates for $\bL^p$ norms and Besov spa\-ces on the sphere}
We describe the performances of the procedure by the following theorem. Remark
that the condition $s>\frac dr$ implies $f\in B^s_{r,q}\subset
B^{s-\frac dr}_{\infty,q}$ so that $f$ is continuous. By $\bE_f$ we denote the expectation
taken with respect to a probability with respect to which the r.v.'s $(X_n)_n$ are
i.i.d. with common density $f$.
\begin{theorem} \label{mM}
For $0< r\le \infty$, $p\ge 1$,  $s>\frac dr$ we have

a) For any $z>1$, there exist some constants $c_\infty=c_\infty(s,p,r,A,M)$ such that if $\kappa>\frac{z+1}6$,
%
\begin{align}
&\sup_{f\in B^s_{r, q}(M)}\bE_f \|\hf -f\|_\infty^z\le
 c_\infty(\log n)^{z-1}\left[\frac n{\log
n}\right]^{\frac{-(s-\frac d r) }{2(s-d(\frac 1r -\frac 12)}}\ .\label{mMinfty}
\end{align}
b) For $1\le p< \infty$ there exist some constant $c_p=c_p(s,r,p,A,M)$ such that if $\kappa>\frac{p}{12}$,
\begin{align}
&\sup_{f\in B^s_{r, q}(M)}\bE_f \|\hf -f\|_p^p\le c_p(\log
n)^{\alpha_p}
\left[\frac n{\log n}\right]^{-\frac{(s-d(\frac1r-\frac 1p))p}
{2(s-d(\frac 1r-\frac12))}},
\label{mMpi}
\end{align}
where $\alpha_p=p-1+1_{\{r=\frac{dp}{2s+d}\}}$, if $r\le\frac{dp}{2s+d}$, whereas
\begin{align}
&\sup_{f\in B^s_{r, q}(M)}\bE_f \|\hf -f\|_p^p\le c_p(\log n)^{p -1}
\left[\frac n{\log n}\right]^{\frac{-sp }{2s+d}},\quad\hbox{if }
r>\frac{dp}{2s+d}\cdotp
\label{mMpidense}
\end{align}
%
%
\end{theorem}
\begin{figure}[h]
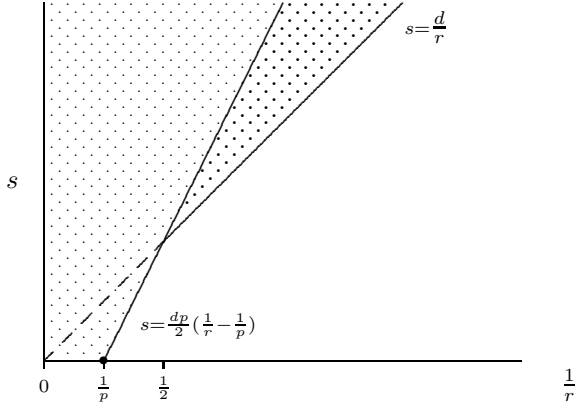

\hbox to\hsize\bgroup\hss
\beginpicture
\setcoordinatesystem units <1.25in,1.25in>
\setplotarea x from 0 to 2, y from 0 to 1.5
\axis bottom  ticks short withvalues $\scriptstyle0$ $\scriptstyle\frac 1p$ $\scriptstyle\frac 12$ / at 0 .25 .5 / /
\axis left label {$s$} shiftedto x=0 
/
\plot .5 .5 1.5 1.5 /
\plot .25 0 1 1.5  /
\setdashes
\plot 0 0 .5 .5 /
\put{${\scriptstyle s=}\frac dr$} at 1.6 1.4
\put{${\scriptstyle s=\frac {dp}2(\frac 1r-\frac 1p)}$} [l] at .4 .15
\setshadesymbol <.1pt,.1pt,.1pt,z> ({\tiny.})
\setshadegrid span <3pt>
\vshade  0 0 1.5 0.25 0 1.5 1 1.5 1.5 /
\setshadesymbol <.1pt,.1pt,.1pt,z> ({.})
\hshade .5 .5 .5 1.5 1 1.5 /
\put{$\scriptstyle\bullet$} at 0.25 0
\put {$\frac1r$} at 2.2 -.1
\endpicture
\hss\egroup
\caption{\footnotesize Type of behaviour of the minimax estimator as a function of $s,r$ and $p$.
The region with the large dots is the one corresponding to the sparse case of
(\ref{mMpi}). The region with the small dots corresponds to the regular case so that (\ref{mMpidense}) holds (see also
Remark \ref{discontinuous}. It should be noticed that
if $p\le 2$ then the slope of the straight line starting at $\frac 1p$ is smaller than the other one, so that the sparse region is empty.\label{sparse-regular}}
\end{figure}

\begin{rem}\rm
Usually the case (\ref{mMpi}) is referred to as the {\it sparse case}, whereas (\ref{mMpidense}) is the {\it regular case}.
Remark that if $p\le 2$, then we are always in the regular case (see also Figure \ref{sparse-regular}).
\end{rem}
\begin{rem}\label{discontinuous}
\rm A closer look to the proof shows that, in the regular case, if we assume
$\|f\|_\infty\le M$, then we can drop the restriction $s>\frac dr$ without
any modification if $1\le p\le 2$. In the case $p>2$,
using an additional modification allowing $J$ to depend also in $p$,
$2^J=\big(\frac n{\log n}\big)^{\frac p{p-2}}$ to be precise, we obtain the same rate
under the same conditions as in the lower bound (up to logarithmic terms).
\end{rem}
\begin{theorem}\label{lb-fin}(Lower bound)
a) If $1\le p\le 2$,
$$
\sup_{f\in B^s_{rq}(M)}\EE_f(\|\hat f-f\|_p^p)\ge c\, n^{-\frac {sp}{2s+d}}\ .
$$
b) If $2<p\le +\infty$
$$
\sup_{f\in B^s_{rq}(M)}\EE_f(\|\hat f-f\|_p^p)
\ge\begin{cases}c\,  n^{-\frac {sp}{2s+d}} &\mbox{if } s>p\frac d2(\frac 1r-\frac 1p)_+\cr
c\,  n^{-\frac {p(s+d(\frac 1p-\frac 1r))}{2(s+d(\frac 12-\frac 1r))}}
&\mbox{if } \frac dr<s\le p\frac d2(\frac 1r-\frac 1p)\ .\cr
\end{cases}
$$
\end{theorem}
\begin{rem}\rm
As already remarked, up to logarithmic terms, the rates observed are minimax. It is
known that in this kind of estimation, full adaptation yields unavoidable extra
logarithmic terms. The rates of the logarithmic terms obtained  in
Theorem \ref{mM} are suboptimal (for instance, for obvious reason the case
$p=2$ yields much less logarithmic terms). We have focused on a simple proof giving all
the results in a rather clear and readable way. However, using a more intricate proof, the rates could be
improved up to be comparable with  those in \cite{MR1394974}.
\end{rem}
\section{Simulations}\label{simulations}
In this section we produce the result of two numerical experiments on the sphere
$\mathbb{S}^2$.
In both of them the major question concerns the choice of the values
of $J$ and $\kappa$. Actually in practical (finite sample)
situations the values given in Theorem \ref{mM}
should be considered just as a reasonable hint. The sets of
cubature points in the simulations that follow have been taken from
the web site of R.~Womersley
\verb@http://web.maths.unsw.edu.au/~rsw@.

\begin{figure}
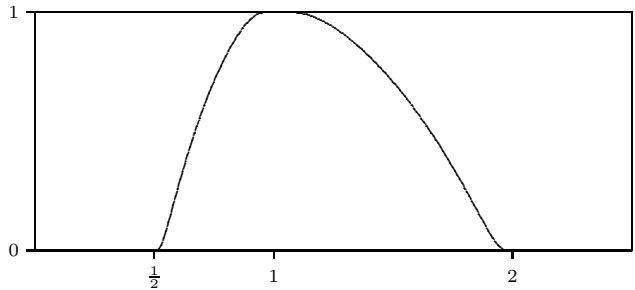

\hbox to\hsize\bgroup\hss
\beginpicture
\setcoordinatesystem units <1.25truein,1.25truein> \setplotarea x
from 0 to 2.5, y from 0 to 1 \axis bottom  ticks short withvalues
{$\scriptstyle\frac 12$} {$\scriptstyle1$}  {$\scriptstyle2$} / at
.5 1 2 / / \axis left shiftedto x=0 ticks short withvalues
{$\scriptstyle 0$} {$\scriptstyle1$} / at 0 1 / / \axis right  /
\axis top / \plot
    0.      0.
    0.01    0.
    0.02    0.
    0.03    0.
    0.04    0.
    0.05    0.
    0.06    0.
    0.07    0.
    0.08    0.
    0.09    0.
    0.1     0.
    0.11    0.
    0.12    0.
    0.13    0.
    0.14    0.
    0.15    0.
    0.16    0.
    0.17    0.
    0.18    0.
    0.19    0.
    0.2     0.
    0.21    0.
    0.22    0.
    0.23    0.
    0.24    0.
    0.25    0.
    0.26    0.
    0.27    0.
    0.28    0.
    0.29    0.
    0.3     0.
    0.31    0.
    0.32    0.
    0.33    0.
    0.34    0.
    0.35    0.
    0.36    0.
    0.37    0.
    0.38    0.
    0.39    0.
    0.4     0.
    0.41    0.
    0.42    0.
    0.43    0.
    0.44    0.
    0.45    0.
    0.46    0.
    0.47    0.
    0.48    0.
    0.49    0.
    0.5     0.
    0.51    0.0001345
    0.52    0.0057761
    0.53    0.0233575
    0.54    0.0501737
    0.55    0.0824075
    0.56    0.1174925
    0.57    0.1538836
    0.58    0.1906656
    0.59    0.2272944
    0.6     0.2634459
    0.61    0.2989279
    0.62    0.3336294
    0.63    0.3674891
    0.64    0.4004771
    0.65    0.4325827
    0.66    0.4638071
    0.67    0.4941586
    0.68    0.5236495
    0.69    0.5522938
    0.7     0.5801056
    0.71    0.6070988
    0.72    0.6332859
    0.73    0.6586776
    0.74    0.6832825
    0.75    0.7071068
    0.76    0.7301541
    0.77    0.7524253
    0.78    0.7739179
    0.79    0.7946263
    0.8     0.8145413
    0.81    0.8336496
    0.82    0.8519338
    0.83    0.8693717
    0.84    0.8859362
    0.85    0.9015943
    0.86    0.9163068
    0.87    0.9300278
    0.88    0.9427043
    0.89    0.9542757
    0.9     0.9646742
    0.91    0.9738261
    0.92    0.9816550
    0.93    0.9880890
    0.94    0.9930738
    0.95    0.9965987
    0.96    0.9987405
    0.97    0.9997272
    0.98    0.9999833
    0.99    1.0000000
    1.      1.
    1.01    1.
    1.02    1.0000000
    1.03    0.9999988
    1.04    0.9999833
    1.05    0.9999136
    1.06    0.9997272
    1.07    0.9993562
    1.08    0.9987405
    1.09    0.9978326
    1.1     0.9965987
    1.11    0.9950167
    1.12    0.9930738
    1.13    0.9907646
    1.14    0.9880890
    1.15    0.9850504
    1.16    0.9816550
    1.17    0.9779107
    1.18    0.9738261
    1.19    0.9694107
    1.2     0.9646742
    1.21    0.9596260
    1.22    0.9542757
    1.23    0.9486322
    1.24    0.9427043
    1.25    0.9365002
    1.26    0.9300278
    1.27    0.9232944
    1.28    0.9163068
    1.29    0.9090714
    1.3     0.9015943
    1.31    0.8938808
    1.32    0.8859362
    1.33    0.8777651
    1.34    0.8693717
    1.35    0.8607601
    1.36    0.8519338
    1.37    0.8428960
    1.38    0.8336496
    1.39    0.8241972
    1.4     0.8145413
    1.41    0.8046837
    1.42    0.7946263
    1.43    0.7843706
    1.44    0.7739179
    1.45    0.7632692
    1.46    0.7524253
    1.47    0.7413868
    1.48    0.7301541
    1.49    0.7187275
    1.5     0.7071068
    1.51    0.6952919
    1.52    0.6832825
    1.53    0.6710779
    1.54    0.6586776
    1.55    0.6460806
    1.56    0.6332859
    1.57    0.6202924
    1.58    0.6070988
    1.59    0.5937037
    1.6     0.5801056
    1.61    0.5663028
    1.62    0.5522938
    1.63    0.5380766
    1.64    0.5236495
    1.65    0.5090108
    1.66    0.4941586
    1.67    0.4790913
    1.68    0.4638071
    1.69    0.4483046
    1.7     0.4325827
    1.71    0.4166403
    1.72    0.4004771
    1.73    0.3840931
    1.74    0.3674891
    1.75    0.3506669
    1.76    0.3336294
    1.77    0.3163809
    1.78    0.2989279
    1.79    0.2812790
    1.8     0.2634459
    1.81    0.2454442
    1.82    0.2272944
    1.83    0.2090233
    1.84    0.1906656
    1.85    0.1722663
    1.86    0.1538836
    1.87    0.1355931
    1.88    0.1174925
    1.89    0.0997089
    1.9     0.0824075
    1.91    0.0658031
    1.92    0.0501737
    1.93    0.0358764
    1.94    0.0233575
    1.95    0.0131448
    1.96    0.0057761
    1.97    0.0015720
    1.98    0.0001345
    1.99    0.0000002
    2.      0.
    2.01    0.
    2.02    0.
    2.03    0.
    2.04    0.
    2.05    0.
    2.06    0.
    2.07    0.
    2.08    0.
    2.09    0.
    2.1     0.
    2.11    0.
    2.12    0.
    2.13    0.
    2.14    0.
    2.15    0.
    2.16    0.
    2.17    0.
    2.18    0.
    2.19    0.
    2.2     0.
    2.21    0.
    2.22    0.
    2.23    0.
    2.24    0.
    2.25    0.
    2.26    0.
    2.27    0.
    2.28    0.
    2.29    0.
    2.3     0.
    2.31    0.
    2.32    0.
    2.33    0.
    2.34    0.
    2.35    0.
    2.36    0.
    2.37    0.
    2.38    0.
    2.39    0.
    2.4     0.
    2.41    0.
    2.42    0.
    2.43    0.
    2.44    0.
    2.45    0.
    2.46    0.
    2.47    0.
    2.48    0.
    2.49    0.
    2.5     0.
/
\endpicture
\hss\egroup
\caption{\footnotesize The function $b$.\label{fig-bb}}
\end{figure}

We realized the function $\varphi$ of \S\ref{ssec-littlewood} by
connecting the levels $0$ and $1$ with a function that is the
primitive, suitably rescaled of the function $x\to
e^{-(1-x^2)^{-1}}$, set to be equal to $0$ outside $[-1,1]$. The
shape of the resulting function $b$ is given in Figure \ref{fig-bb}.
For this choice of $b$, we have
$$
\frac18\int_{1/2}^2 tb^2(t)\, dt \simeq 0.107
$$
which, as remarked above gives an indication about the square of the
value of the $L^2$ norm of a needlet
$\psi_{j\xi}$.
In both the examples below we considered samples of cardinality
$n=2000$ and $n=8000$. The hint for the value of $J$ of Theorem
\ref{mM} is $J=\frac 12\,\log_2\big(\frac n{\log n}\big)$, which
gives the values $J\sim 4.02$ and $J\sim4.9$ respectively. One
should keep in mind that at a given level $j$ it is necessary to have enough
cubature points in order to integrate exactly all polynomials up to
the degree $2(2^{j+1}-1)=2^{j+2}-2$, which means $\sim2^{2j+4}$
cubature points with  Womersley's set (recall that on the sphere the
polynomials of degree  $d$ form a vector space of dimension
$(2d+1)^2$). This gives $2^{10}=1024$ cubature points for $j=3$, $2^{12}=4096$ for
$j=4$ and $2^{14}=16384$ for $j=5$. To avoid to have
more coefficients than observations, we
decided to set $J=3$ for $n=2000$ and $J=4$ for $n=8000$.

As for the value of $\kappa$, we shall give the result with
$\kappa=k_0\sqrt{0.107}\,M$, where $M$ is an a bound for
$\|f\|_\infty$, trying different values of $k_0$. Recall that
this means that the threshold kills all coefficients $\beta_{j\xi}$
such that $|\beta_{j\xi}|<\kappa\sqrt{\frac {\log n}n}$
\begin{example}\rm
$f=\frac 1{4\pi}$, the uniform density. In this case in the development
\tref{develop} it holds $\beta_{j\xi}=\langle f,\psi_{j\xi}\rangle_{L^2}=0$
for every $j$ and $\xi$. Therefore a first simple way of assessing the performance of the procedure
is to count the number of coefficients that survive thresholding.
Of course in this case a good
estimate is such that the coefficients $\beta_{j,\xi}$
fall below the threshold. Taking into account
Lemma \ref{lemma-inequality} the square root of the sum of the squares
of the coefficients surviving thresholding gives an estimate of $\|\hat f-f\|_2$.
Therefore a measure of the goodness of the fit is obtained by
taking the sum of their squares.
Tables \ref{tab1} and \ref{tab2} give the number of surviving coefficients for
different values of the constant $k_0$.
\begin{table}
\begin{center}
\begin{tabular}{|l|r|c|c|c|}
  \hline
 & $j=0$ & $j=1$& $j=2$& $j=3$\\  \hline
  $k_0=1$ & 8 (.89) & 29 (.45)& $96$ ($.38$)&$471$ ($.46$) \\
  \hline
  $k_0=1.5$&7 (.78) & 16 (.25)& $45$ ($.18$)&$264$ ($.26$) \\
    \hline
  $k_0=2$&$4$ ($.44$) & $4$ ($.06$)& $29$ ($.11$)&$126$ ($.12$) \\
  \hline
\end{tabular}
\end{center}
\caption{number of coefficients surviving thresholding for various values of $k_0$, $n=2000$.\label{tab1}}
\end{table}
\begin{table}
\begin{center}
\begin{tabular}{|l|r|c|c|c|c|}
  \hline
 & $j=0$ & $j=1$& $j=2$& $j=3$& $j=4$\\  \hline
  $k_0=1$ & 4 (.44) & 28 (.44)& $96$ ($.38$)&$413$ ($.40$)&$1610$ ($.39$) \\
  \hline
  $k_0=1.5$&2 (.22) & 12 (.19)& $50$ ($.20$)&$207$ ($.20$)&$921$ ($.20$) \\
    \hline
  $k_0=2$&$1$ ($.11$) & $4$ ($.06$)& $16$ ($.06$)&$97$ ($.09$)&$368$ ($.09$) \\
  \hline
\end{tabular}
\end{center}
\caption{number of coefficients surviving thresholding for various values of $k_0$, $n=8000$.\label{tab2}}
\end{table}
%
In order to kill all the coefficients one should choose $k_0=5.4$ for $n=2000$ and $k_0=2.8$
for $n=8000$.
The estimate of the $L^2$ norm of the difference between $\hat f$
and $f$ by taking the square root of the sum of the squares of the coefficients is
\begin{center}
\begin{tabular}{|c|c|c|c|}
\hline
$k_0$&1&1.5&2\\
\hline
$n=2000$&0.146&0.131&0.107\\
\hline
$n=8000$&0.108&0.0834&0.060\\
\hline
\end{tabular}
\end{center}
\end{example}
\begin{example}\rm Let us consider a mixture $f$ of two densities
of the form $f_i(x)=c_i \exp(-k_i|x-x_i|^2)$, $i=1,2$, for $k_1=.7$ and
$k_2=2$ and with weights $0.65$ and $0.35$ respectively. Here the centers $x_i$ of the
two bell-shaped densities were taken to be $x_1=(0,1,0)$, $x_2=(0,-.8,.6)$. With these choices
it turns out that $\|f\|_\infty=0.26$. The graph
of $f$ in the coordinates $(\phi,\theta)$ ($\phi=$longitude, $\theta=$colatitude)
is given in Figure \ref{two-bumps}.

\begin{figure}[h]
\centering
\includegraphics[width=11cm]{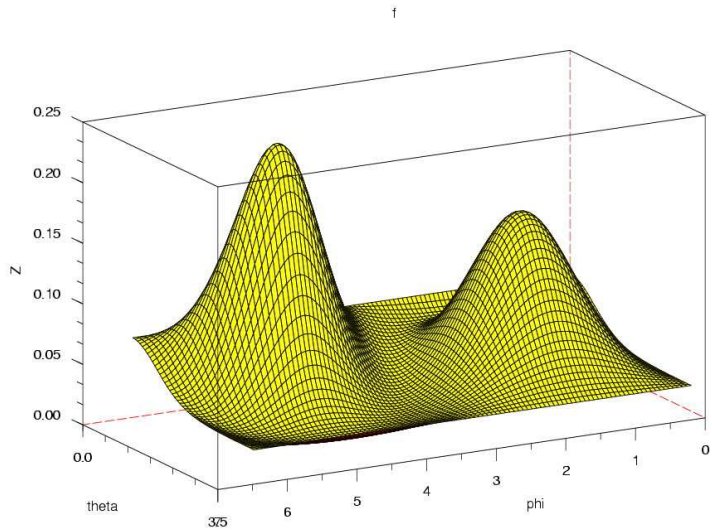}
\caption{Graph of the density with two bumps.\label{two-bumps}}
\end{figure}
The estimator $\hat f$ obtained with the choice $k_0=1.1$ has the graph of Figure \ref{two-bumps-hat32000}.
\begin{figure}[h!]
\centering
\includegraphics[width=11cm]{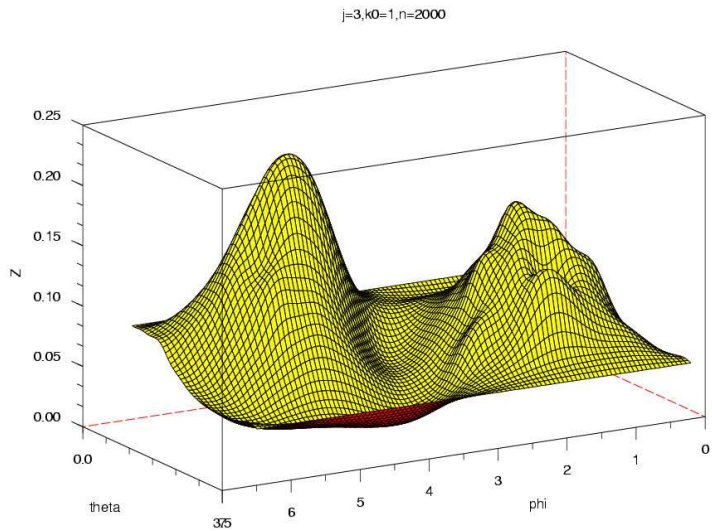}
\caption{Graph of the estimated density ($n=2000$, $k_0=1.1$).\label{two-bumps-hat32000}}
\end{figure}
If one chooses $k_0=1.5$ the graph becomes the one of Figure \ref{two-bumps-hat22000}.
At a closer inspection it turns out that with this value of $k$ all coefficients at level $j=3$ do not pass the thresholding.
It looks very much like the graph of $f$, even though some differences in shape are apparent.
An estimate of the norm $\|\hat f-f\|_\infty$ computed on a grid gives $\|\hat f-f\|_\infty\sim 0.054$.
\begin{figure}[h!]
\centering
\includegraphics[width=11cm]{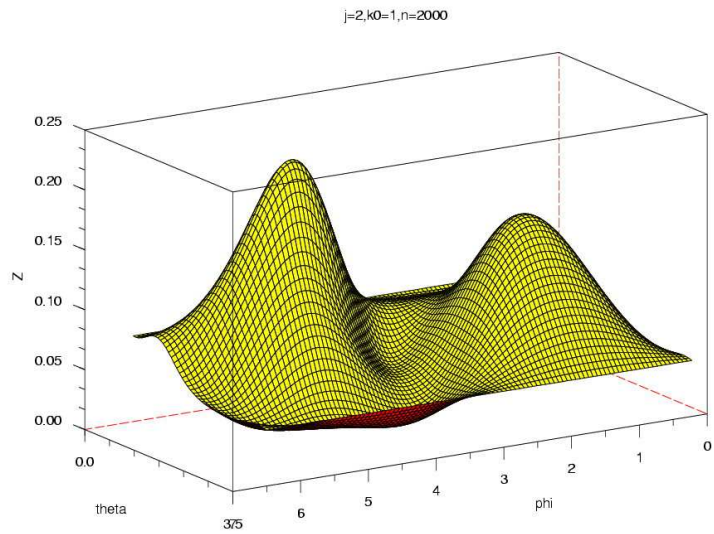}
\caption{Graph of the estimated density ($n=2000$, $k_0=1.65$).\label{two-bumps-hat22000}}
\end{figure}
We repeated the simulation with $n=8000$ observations. The results are reported in Figures \ref{two-bumps-hat48000-11}
and \ref{two-bumps-hat48000-165} and are to be considered rather satisfactory. It should be stressed that a very limited
number of coefficients passes thresholding at a frequency $j>2$. This behaviour is expected: $f$ being very regular, it belongs
to a space $B^{\infty}_{\infty,\infty}$ and thus its needlet coefficients decay very rapidly.
\begin{figure}[h!]
\centering
\includegraphics[width=11cm]{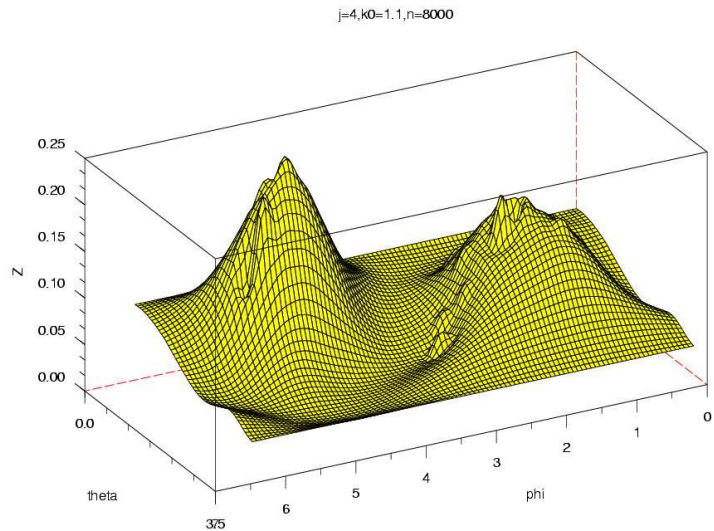}
\caption{Graph of the estimated density ($n=8000$, $k_0=1.1$).\label{two-bumps-hat48000-11}}
\end{figure}
\begin{figure}[h!]
\centering
\includegraphics[width=11cm]{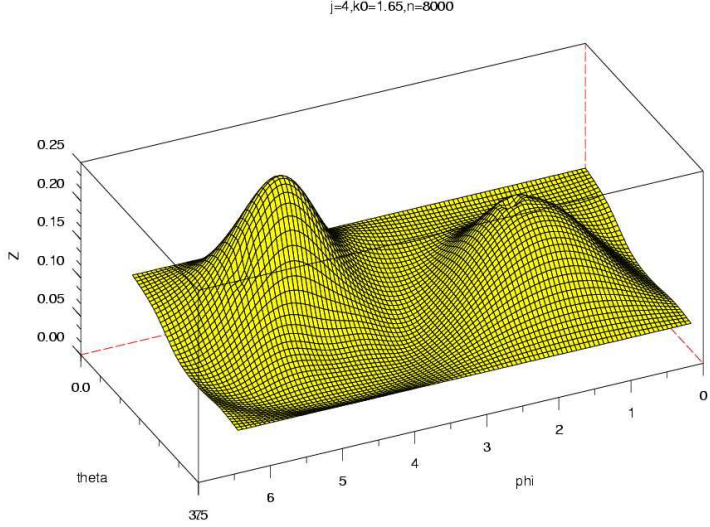}
\caption{Graph of the estimated density ($n=8000$, $k_0=1.65$). It is rather satisfactory, but
for a small dent on the top of the lowest bump. With this value of $k_0$ only 2 coefficients pass thresholding for $j=4$
and no one at level $j=3$. Now
$\|\hat f-f\|_\infty\sim 0.028$\label{two-bumps-hat48000-165}}.
\end{figure}
\end{example}%
\section{Proof of Theorem \protect{\ref{mM}}}
In the sequel we note $t(\hat \beta_{j,\xi})=
\hbjk\,{1}_{\{|\hbjk|\geq \,\kappa\,c_n\}}$, so that the needlet estimator
\tref{eq:need} is
$$
\hf=\frac 1{|\mathbb S^d|}+\sum_{j=0}^J\sum_{\eta\in
\ccl Z_j} \,t(\hbjk)\,\psijk\ .
$$
In this section and in the next one the density $f$ is fixed and
we shall write $\bE$ instead of $\bE_f$,
as there is no danger of confusion.

The following proposition collects the main estimates needed in the proof.
\begin{proposition}\label{key}
Let $J_1\le J$ be such that,  for all $J_1\le j\le J$, $|\bjk|\le \tfrac\kappa 2 t_n$
(possibly $J_1=J$; obviously, when $f$ belongs to a Besov class,
$J_1$ depends on the ``regularity'' $s$).
Then for any $\gamma>0, s>0, z\ge 1$, we have

1) if $\kappa>\frac13(\frac\gamma d+\frac12)$
\begin{equation}\label{propsup}
\begin{array}{c}\displaystyle\sum_{j=0}^{{J}}2^{\gamma j}\bE \Bigl[\sup_\eta|t(\hat\beta_{j\eta})-\bjk|^z\Bigr]
\le\\
\displaystyle\le C\Bigl[2^{J_1\gamma}(J_1+1)^zn^{-\frac z2}+
\sum_{j=J_1+1}^{{J}}2^{\gamma j}\sup_{\eta\in \Z_j}  |\bjk|^z
+n^{-\frac z2}\Big]\ .\end{array}
\end{equation}
2) if $\kappa>\frac\gamma {6d}$
\begin{equation}\label{propsum1}
\begin{array}{c}
\displaystyle\sum_{j=0}^{{J}}2^{(\gamma-d) j}\bE \sum_\eta|t(\hat\beta_{j\eta})-\bjk|^z
\le\\
\displaystyle\le  C\sum_{j=0}^{{J_1}}2^{(\gamma-d) j}n^{-\frac z2}\sum_{\eta\in \Z_j}1_{\{|\bjk|>\frac\kappa 2t_n\}}+\qquad\qquad\qquad\\
\displaystyle\qquad\qquad\qquad+C\sum_{j=0}^{{J_1}}2^{(\gamma-d) j}
\sum_{\eta\in \Z_j}1_{\{|\bjk|\le 2\kappa t_n\}}+Cn^{-\frac z2}
\end{array}
\end{equation}
and 
\begin{equation}\label{propsum}
\begin{array}{c}
\displaystyle\sum_{j=0}^{{J}}2^{(\gamma-d) j}\bE \sum_\eta|t(\hat\bjk)-\bjk|^z
\le\\
\displaystyle\le C \Big[2^{J_1\gamma}(J_1+1)^zn^{-\frac z2}+\!\!
\sum_{j=J_1+1}^{{J}}2^{(\gamma-d) j}\!\!\sum_{\eta\in \Z_j}  |\bjk|^z
+n^{-\frac z2}\Big]\ .
\end{array}
\end{equation}
\end{proposition}

We delay the proof of Proposition \ref{key} to \S\ref{sect-key} and derive from it the
proof of Theorem \ref{mM}.
In this proof $C$ will denote an absolute constant which may change from line to line.
Let us now prove that Proposition \ref{key} yields to the statements of Theorem \ref{mM}:

Let us prove the $L^\infty$ upper bound (\ref{mMinfty}), first under the condition
$q=r=\infty$
\begin{equation}\label{first}
\begin{array}{c}
\displaystyle\EE\|\hat f-f\|_\infty^z \le\\
\displaystyle C\Big[\bE\Big\|\sum_{j=0}^{{J}}\sum_{\eta\in \Z_j}
(t(\hat\beta_{j\xi})-\bjk)\psijk\Big\|_\infty^z +\Big\|\sum_{j>{J}}\sum_{\eta\in
\Z_j}  \bjk\psijk\Big\|_\infty^z\Big]\\
:= I+II\ .
\end{array}
\end{equation}
The term $II$ is easy to analyze:
as $f$ belongs to $B^s_{\infty,\infty}(M)$, we have using  (\ref{lp}) and \ref{inequality},
$$
\displaylines{
\Big\|\sum_{j>{J}}\sum_{\eta\in \Z_j}  \bjk\psijk\Big\|_\infty\le
\sum_{j>{J}}\Big\|\sum_{\eta\in \Z_j}  \bjk\psijk\Big\|_\infty\le C
\sum_{j>{J}}\sup_{\eta\in \Z_j}  |\bjk|\|\psijk\|_\infty\le\cr
\le C
\sum_{j>{J}} 2^{-j[s+\frac d{2}]}2^{\frac{jd}2}\le \left(\frac{\log
n}n\right)^\frac sd\ .\cr
}
$$
Then we only need to remark that $\frac{s}{d}\ge\frac{s}{d+2s} $ for $s>0$.


As for $I$, using the triangular inequality together with H\"older inequality,
then (\ref{lp}), and  (\ref{propsup}) with $\gamma =\frac {dz}2, z>1$,
 we get
$$
\displaylines{
I\le C(J+1)^{z-1}\sum_{j=0}^{{J}}2^{\frac{jdz}2}\bE\sup_{\eta\in \Z_j}
|t(\hat\bjk)-\bjk|^z\le\cr
\le
C (J+1)^{z-1}[2^{J_1\frac {dz}2}(J_1+1)^zn^{-\frac z2}+C
\sum_{j=J_1+1}^{{J}}2^{\frac {dz}2 j}\sup_{\eta\in \Z_j}  |\bjk|^z
+n^{-\frac z2}]\ .\cr
}
$$
As $f$ belongs to  $B^s_{\infty,\infty}(M)$,  $|\bjk|\le M2^{-j(s+\frac d2)}$ and
we can see that
$$
2^{J_1}=
\frac \kappa{2M}
\Big[\frac n {\log n}\Big]^{\frac 1{2s+ d}},
$$
is adequate, it is easy to conclude.

For arbitrary $q, r $  (\ref{mMinfty}) is now easy to deduce from the
previous computation by the Besov  embedding (Theorem \ref{Besov-embedding})
$
B^s_{r,q}(M)\subset B^{s-\frac dr}_{\infty,\infty}(M).
$
Let us prove (\ref{mMpidense}), that is the regular case.
We observe first that since $B^s_{r,q}(M)\subset B^s_{p,q}(M)$ for $r\ge p$,
this case will be assimilated to the case $p=r$, and from now on, we only consider
$r\le p $.
We follow the same arguments as above. (\ref{first}) can be replaced by
\begin{equation}\label{first'}
\begin{array}{c}
\displaystyle\EE\|\hat f-f\|_p^p \le\\
\displaystyle\le C\Bigl(\bE\Big\|\sum_{j=0}^{{J}}\sum_{\eta\in \Z_j}
(t(\hat\bjk)-\bjk)\psijk\Big\|_p^p +\Big\|\sum_{j>{J}}\sum_{\eta\in \Z_j}
\bjk\psijk\Big\|_p^p\Bigr)=:\\
\displaystyle=: I+II\ .
\end{array}
\end{equation}
For $II$ using the embedding $B^s_{r,q}(M)\subset B^{s-\frac dp+\frac dp}_{p,q}(M)$,
for $r\le p $, we have
$$
II^{\frac 1p}\le C\Big\|\sum_{j>{J}} \sum_{\eta\in \Z_j}  \bjk\psijk\Big\|_p
\le C2^{-J(s-\frac dr+\frac dp)}.
$$
And it is easy to verify that $\frac sd -\frac 1r+\frac 1p\ge \frac s{2s+d}$ on the
zone that we are considering in this part.
In effect as $s\ge \frac{p }2(\frac dr-\frac dp)$,
$\frac s{2s+d}\le \frac{sr}{dp}$ we have
$\frac sd-\frac 1r +\frac 1p-\frac{sr}{dp}=(\frac 1r-\frac 1p)(\frac sd r-1)\ge 0$.

For $I$, we have using the triangular inequality together with H\"older inequality,
$$
\displaylines{
\bE\Big\|\sum_{j=0}^{{J}}\sum_{\eta\in \Z_j}  (t(\hat\bjk)-\bjk)\psijk\Big\|_p^p\le\cr
\le
C(J+1)^{p-1}\sum_{j=0}^{{J}}2^{jd(\frac p 2 -1)}\sum_{\eta\in \Z_j}
\bE|t(\hat\bjk)-\bjk)|^p.\cr
}
$$
Then we need only to use (\ref{propsum}), with $\gamma=\frac {dp} 2 , z=p$, to obtain
$$
I\le C(J+1)^{p-1} \Big[2^{J_1d\frac p 2}(J_1+1)^p n^{-\frac p 2}+
\sum_{j=J_1+1}^{{J}}2^{d(\frac p 2-1) j}\sum_{\eta\in \Z_j}  |\bjk|^p
+n^{-\frac p2}\Big].
$$
It is easy to realize that again
$$2^{J_1}=
\frac \kappa{2M}
\Big[\frac n {\log n}\Big]^{\frac 1{2s+ d}},$$
is adequate, and to observe that the first term in the sum has the right order.
For the second term, it can be bounded (as $p\ge r$) by
$$
\displaylines{
C(J+1)^{p-1}\sum_{j=J_1+1}^{{J}}2^{d(\frac p 2-1) j}
\sum_{\eta\in \Z_j}  |\bjk|^r[\kappa t_n]^{p - r}\le\cr
\le C(J+1)^{p-1}\sum_{j=J_1+1}^{{J}}2^{d(\frac p 2-1) j}M^p
2^{-jr(s+\frac d2-\frac dr)}t_n^{p - r}\le\cr
\le Ct_n^{p - r}
2^{-J_1(sr-\frac d2(p -r))}\cr
}
$$
as, $sr-\frac d2(p -r)\ge 0.$
Now, this term obviously is of the right order.

Again we proceed as above and observe first that in order to have $s>0$ as well as $s\le \frac{p d}2(\frac 1r-\frac 1p)$, it is necessary that $p\ge r$:
\begin{equation}\label{first''}
\begin{array}{c}
\displaystyle\EE\|\hat f-f\|_p^p \le\\
\displaystyle\le C\bE\|\sum_{j=0}^{{J}}\sum_{\eta\in \Z_j}
(t(\hat\bjk)-\bjk)\psijk\|_p^p +\|\sum_{j>{J}}\sum_{\eta\in \Z_j}  \bjk\psijk\|_p^p:=\\
\displaystyle=: I+II\ .
\end{array}
\end{equation}
For $II$ using the embedding, $B^s_{r,q}(M)\subset B^{s-\frac dr+\frac dp}_{p,q}(M)$, for $r\le p $, we have:
$$
II^{\frac 1p}\le C\Big\|\sum_{j>{J}} \sum_{\eta\in \Z_j}  \bjk\psijk\Big\|_r\le C2^{-J(s-\frac dr+\frac dp)}.
$$
And it is easy to verify that $\frac sd -\frac 1r+\frac 1p\ge \frac{(s-d(\frac1r-\frac 1p))}{2(s-d(\frac 1r-\frac12))}$,
since  $2(s-d(\frac 1r-\frac12))\ge d$, when $s>\frac dr$.

For $I$, again, we have using the triangular inequality together
with H\"older inequality, 
$$
\displaylines{
\bE\Big\|\sum_{j=0}^{{J}}\sum_{\eta\in \Z_j}  (t(\hat\bjk)-\bjk)\psijk\Big\|_p^p\le \cr
\le C(J+1)^{p-1}\sum_{j=0}^{{J}}2^{jd(\frac p 2 -1)}\sum_{\eta\in \Z_j}
\bE|t(\hat\bjk)-\bjk)|^p.\cr
}
$$
Then we need  to use (\ref{propsum1}), with $\gamma=d\frac p 2 ,z=p$, to obtain:
\begin{align*}
I&\le C(J+1)^{p-1} \sum_{j\le J_1}2^{jd(\frac p 2-1)}\sum_{\eta\in \Z_j} 1_{\{|\bjk|>\frac \kappa 2
t_n\}}|\bjk|^r[\frac \kappa 2
t_n]^{-r}n^{-\frac p 2}+\\
&\qquad\qquad+ C(J+1)^{p-1} \sum_{j\le J}2^{jd(\frac p 2-1)}\sum_{\eta\in \Z_j} 1_{\{|\bjk|\le 2 \kappa
t_n\}}|\bjk|^p+n^{-z/2}\le
\\
&\le 2C(J+1)^{p-1} \sum_{j\le J_1}
2^{jd(\frac p 2-1)}M^r
2^{-jr\{s+d(\frac 12-\frac 1r)\}}n^{\frac{r-p}2}+
\\
&\qquad\qquad+ C(J+1)^{p-1} \sum_{j\ge J_1}2^{jd(\frac p 2-1)}
\sum_{\eta\in \Z_j} 1_{\{|\bjk|\le 2 \kappa
t_n\}}|\bjk|^p +n^{-z/2}\le \\
&\le 2C(J+1)^{p-1}2^{J_1\{d(\frac p2-\frac r2)-sr\}}n^{\frac{r-p}2}+
\\
&\qquad\qquad+ C(J+1)^{p-1} \sum_{j\ge J_1}2^{jd(\frac p 2-1)}
\sum_{\eta\in \Z_j} 1_{\{|\bjk|\le 2 \kappa
t_n\}}|\bjk|^p+n^{-z/2}
\end{align*}
as we are in the sparse region.
It is easy to realize that now, again because we are in the sparse region
$$2^{J_1}=
\frac \kappa{2M}
[\frac n {\log n}]^{\frac 1{2s+ 2d(\frac 12-\frac1r)}}\ ,
$$
is adequate, and to observe then that the first term in the sum has the right order.
For the second term,
let us introduce
$$m:= \frac{d(\frac p2-1)}{s+\frac d2-\frac dr}.$$
We easily observe that $p-m= \frac{p(s+\frac dp-\frac dr)}{s+\frac d2-\frac dr}>0$, and that $m-r=\frac{-sr+(p-r)\frac d2}{s+\frac d2-\frac dr}\ge 0.$
Then, as $B^s_{r,q}(M)\subset B^{s-\frac dr+\frac dm}_{m,q}(M)$
$$
\displaylines{
\sum_{J\ge j\ge J_1}2^{jd(\frac p 2-1)}\sum_{\eta\in \Z_j} 1_{\{|\bjk|\le 2 \kappa
t_n\}}|\bjk|^p \le \cr
\le\sum_{J\ge j\ge J_1}2^{jd(\frac p 2-1)}\sum_{\eta\in \Z_j}|\bjk|^mt_n^{p-m}\le
\sum_{J\ge j\ge J_1}M^m t_n^{p-m}\le\cr
\le J t_n^{\frac{p(s+\frac dp-\frac dr)}{s+\frac d2-\frac dr}},\cr
}
$$
which gives the right order. Observe that the term $J$ (which is of logarithmic order), can be avoided by choosing $\tilde m$ instead of $m$ in such a way that $\tilde m>m$, but $r<\tilde m$. This can be done except for the case where $r=\frac{dp}{2s+d}$ where this logarithmic term is unavoidable.
\section{Proof of Proposition \ref{key}}\label{sect-key}

The proof of  Proposition \ref{key} relies on  the following lemma:


\begin{lemma}\label{quattordici}
There exist constants $\sigma^2>0, \; C,\; c$, such that, as soon as $2^{j}\le [\frac n{\log n}]^{\frac 1d}$,
\begin{align}\label{bern}
\bP\{|\hat\bjk-\bjk|\ge v\}&\le
2\exp\Big\{-\frac{nv^2}{2(\sigma^2 +\frac13 {vc2^{\frac{jd}2}})}\Big\},\quad \forall \; v>0,
\\
\bE|\hat\bjk-\bjk|^q &\le  s_qn^{-\frac q2},\quad \forall \; q\ge 1\label{mom}\\
\bE\sup_{\eta}|\hat\bjk-\bjk|^q &\le  s_q'(j+1)^qn^{-\frac q2},\quad \forall \; q\ge 1\label{momsup}\\
\bP(|\hat\bjk-\bjk|\ge \tfrac{\kappa} 2  t_n{})&\le C
2n^{-6\kappa},\quad \forall \; \kappa\ge 6\sigma^2\label{gd}
\end{align}
\end{lemma}
\noindent{\bf Proof of the lemma}
(\ref{bern}) is simply Bernstein inequality, noticing that
$$
\bE(\psijk(X_i))^2\le \|f\|_\infty\|\psijk\|_2^2\le MC =:\sigma^2
$$
and $\|\psijk(X_i)\|_\infty\le c2^{\frac{jd}2}$.
The following inequality directly follows  from (\ref{bern}), when $2^{jd}\le n$.
\begin{equation}\label{bern1}
\bP\{|\hat\bjk-\bjk|\ge v\}\le 2[e^{-\frac{nv^2}{4\sigma^2}}
+e^{-\frac{3\sqrt{n}v}{ 4c}}
]
\end{equation}
(\ref{mom}) follows from (\ref{bern1}):
$$
\displaylines{
\bE|\hat\bjk-\bjk|^q =\int_{\bR^+_*} v^{q-1}\bP(|\hat\bjk-\bjk|\ge v)\,dv\le\cr
\le \int_{\bR^+_*} v^{q-1}2\Big[e^{-\frac{nv^2}{4\sigma^2}}
+e^{-\frac{3\sqrt{n}v}{ 4c}}
\Big]\,dv\le: s_qn^{-\frac q2},\cr
}
$$
using the change of variables $u=\sqrt{n}v$.

(\ref{momsup}) also follows from (\ref{bern1}): take $a=\max\{\frac {8cd}{3}, {2\sqrt 2 d\sigma}\}$
$$
\displaylines{
\bE\sup_\eta|\hat\bjk-\bjk|^q =
\int_{\bR^+_*} v^{q-1}\bP(\sup_\eta|\hat\bjk-\bjk|\ge v)\,dv\cr
\le\int_{0\le v\le \frac {aj}{\sqrt{n}}} v^{q-1}dv+
2c\int_{v\ge \frac {aj}{\sqrt{n}}} v^{q-1}2^{jd}
\Big[e^{-\frac{nv^2}{4\sigma^2}}+e^{-\frac{3\sqrt{n}v}{ 4c}}
\Big]\,dv\cr
}
$$
Now, if $v\ge \frac {aj}{\sqrt{n}}$, $2^{jd}e^{-\frac{nv^2}{4\sigma^2}}\le
e^{-\frac{nv^4}{8\sigma^2}-\frac{nv^4}{8\sigma^2}+jd}\le e^{-\frac{nv^4}{8\sigma^2}}$.
Similarly $2^{jd}e^{-\frac{3\sqrt{n}v}{ 4c}}\le e^{-\frac{3\sqrt{n}v}{ 8c}}$, so that
$$
\bE\sup_\eta|\hat\bjk-\bjk|^q\le \frac 1q\,\Big[\frac {aj}{\sqrt{n}}\Big]^q+
2c\int_{v\ge \frac {aj}{\sqrt{n}}} v^{q-1}\big[e^{-\frac{nv^4}{8\sigma^2}}+e^{-\frac{3\sqrt{n}v}{ {8c}}}
\big]\,dv
$$\qed
Let us now turn to the proof of the Proposition.
We partition our sum in four regions:
\[
\begin{array}{c}
\displaystyle\sum_{j=0}^{{J}}2^{j\gamma}\bE\sup_{\eta\in \Z_j}  |t(\hat\bjk)-\bjk|^z=\\
=\displaystyle\sum_{j=0}^{{J}}2^{j\gamma}\bE\sup_{\eta\in \Z_j}  |t(\hat\bjk)-\bjk|^z
\big\{1_{\{|\hat\bjk|\ge {\kappa}  t_n{}\}}
+1_{\{|\hat\bjk|< {\kappa}  t_n{}\}}\big\}\le\\
\displaystyle\le  \sum_{j=0}^{{J}}2^{j\gamma}\bE\sup_{\eta\in \Z_j}  |\hat\bjk-\bjk|^z
1_{\{|\hat\bjk|\ge {\kappa}  t_n{}\}}
1_{\{|\bjk|\ge \frac\kappa 2  t_n{}\}}\qquad\qquad\qquad
\\
\displaystyle\qquad\qquad\qquad+  \sum_{j=0}^{{J}}2^{j\gamma}\bE\sup_{\eta\in \Z_j}  |\hat\bjk-\bjk|^z
1_{\{|\hat\bjk|\ge {\kappa}  t_n{}\}}
1_{\{|\bjk|< \frac\kappa 2  t_n{}\}}+
\\
\displaystyle+  \sum_{j=0}^{{J}}2^{{j\gamma}}\bE\sup_{\eta\in \Z_j}  |\bjk|^z
1_{\{|\hat\bjk|< {\kappa}  t_n{}\}1_{\{|\bjk|\ge 2{\kappa}  t_n{}\}}}+\qquad\qquad\qquad
\\
\displaystyle\qquad\qquad\qquad+  \sum_{j=0}^{{J}}2^{{j\gamma}}\bE\sup_{\eta\in \Z_j}  |\bjk|^z
1_{\{|\hat\bjk|< {\kappa}  t_n{}\}}1_{\{|\bjk|< 2{\kappa}  t_n{}\}}=
\\
= :Bb+Bs+Sb+Ss \ .
\end{array}\]

We use extensively Lemma \ref{quattordici} in order to bound separately each of the four terms $Bb, Ss, Sb, Bs$.

Using (\ref{momsup})
\begin{align*}
Bb&\le\sum_{j=0}^{{J}}2^{j\gamma}\bE\sup_{\eta\in \Z_j}  |\hat\bjk-\bjk|^z
 1_{\{|\bjk|\ge \frac\kappa 2  t_n{}\}}\\
&\le
\sum_{j=0}^{{J}}1_{\{\exists\; \eta \in \Z_j,\; |\bjk|\ge \frac\kappa 2    t_n{}\}}2^{j\gamma}\bE\sup_{\eta\in \Z_j}  |\hat\bjk-\bjk|^z\\
&\le
\sum_{j=0}^{{J}}1_{\{\exists\; \eta \in \Z_j,\; |\bjk|\ge \frac\kappa 2    t_n{}\}}2^{j\gamma}s_1'(j+1)^zn^{-\frac z2}
\\
&\le C2^{J_1\gamma}(J_1+1)^zn^{-\frac z2}
\end{align*}
where $J_1$ is chosen such that
for $j\ge J_1$, $|\bjk|\le \frac\kappa 2 t_n$. Also
\begin{align*}
Ss&\le\sum_{j=0}^{{J}}2^{j\gamma}\sup_{\eta\in \Z_j}  |\bjk|^z
1_{\{|\bjk|< 2{\kappa}  t_n{}\}}
\\
&\le\sum_{j=0}^{{J_1}}2^{j\gamma}[2{\kappa}  t_n]^z
+\sum_{j=J_1+1}^{{J}}2^{j\gamma}\sup_{\eta\in \Z_j}  |\bjk|^z
%
&
\end{align*}
which gives the proper rate of convergence. Moreover,
using (\ref{momsup}) and (\ref{gd}),
\begin{align*}
Bs&\le \sum_{j=0}^{{J}}2^{j\gamma}\bE\sup_{\eta\in \Z_j}  |\hat\bjk-\bjk|^z
1_{\{|\hat\bjk-\bjk|\ge \frac\kappa 2   t_n{}\}}
1_{\{|\bjk|< \frac\kappa 2  t_n{}\}}\\
&\le \sum_{j=0}^{{J}}2^{j\gamma}\bE\sup_{\eta\in \Z_j}  |\hat\bjk-\bjk|^z1_{\{\exists\; \eta \in \Z_j,\;|\hat\bjk-\bjk|\ge \frac\kappa 2   t_n{}\}}\\
&\le \sum_{j=0}^{{J}}2^{j\gamma}[\bE\sup_{\eta\in \Z_j}  |\hat\bjk-\bjk|^{2z}]^{\frac12}\bP\{\exists\; \eta \in \Z_j,\;|\hat\bjk-\bjk|\ge \frac\kappa 2   t_n{}\}^{\frac12}\\
&\le \sum_{j=0}^{{J}}2^{j\gamma} [s_2'(j+1)^{2z}n^{-z}]^{\frac12}[c2^{jd}n^{-6\kappa}]^{\frac12}\le n^{-\frac z2}
\end{align*}
where $\kappa>\frac 13(\frac\gamma d+\frac12)$.
Finally, using
(\ref{gd}), and the fact that for $f$ bounded, $|\bjk|\le C2^{-\frac{jd}2}$
\begin{align*}
Sb&\le \sum_{j=0}^{{J}}2^{j\gamma}\bE\sup_{\eta\in \Z_j}  |\bjk|^z
1_{\{|\bjk-\hat\bjk|\ge {\kappa}  t_n{}\}}1_{\{|\bjk|\ge 2{\kappa}  t_n{}\}}
\\
&\le \sum_{j=0}^{{J}}2^{j\gamma}M2^{-jz\frac d2}\bP\{\exists\; \eta \in \Z_j,\;|\bjk-\hat\bjk|\ge {\kappa}  t_n{}\}
\\
&\le \sum_{j=0}^{{J}}[c2^{j[d(1-\frac z2)+\gamma]}n^{-6\kappa}]
\\
&\le C [2^{J[d(1-\frac z2)+\gamma]}n^{-6\kappa}]\le n^{-\frac z2}
,\end{align*}
for $\kappa>\frac 16(\frac\gamma d+1)$.
%
\subsection{Proof of \protect(\ref{propsum1}) and (\ref{propsum})}
This proof follows along the lines of the previous one. (\ref{propsum}) is a consequence of (\ref{propsum1}), and the two inequalities will be proved together.
We again separate the four cases.
$$
\displaylines{ \sum_{j=0}^{{J}}2^{j(\gamma-d)}\bE\sum_{\eta\in \Z_j}
|t(\hat\bjk)-\bjk|^z\le\cr \le
\sum_{j=0}^{{J}}2^{j(\gamma-d)}\bE\sum_{\eta\in \Z_j}
|\hat\bjk-\bjk|^z 1_{\{|\hat\bjk|\ge {\kappa}  t_n\}} 1_{\{|\bjk|\ge
\frac\kappa 2  t_n{}\}}+\cr
\qquad\qquad\enspace\enspace\>+
\sum_{j=0}^{{J}}2^{j(\gamma-d)}\bE\sum_{\eta\in \Z_j}
|\hat\bjk-\bjk|^z 1_{\{|\hat\bjk|\ge {\kappa}  t_n{}\}} 1_{\{|\bjk|<
\frac\kappa 2  t_n{}\}}+\cr
\qquad\quad+  \sum_{j=0}^{{J}}2^{j(\gamma-d)}\bE\sum_{\eta\in \Z_j}  |\bjk|^z
1_{\{|\hat\bjk|< {\kappa}  t_n{}\}}1_{\{|\bjk|\ge 2{\kappa}  t_n{}\}}+\cr
\qquad\qquad+  \sum_{j=0}^{{J}}2^{j(\gamma-d)}\bE\sum_{\eta\in \Z_j}  |\bjk|^z
1_{\{|\hat\bjk|< {\kappa}  t_n{}\}}1_{\{|\bjk|< 2{\kappa}  t_n{}\}}:=\cr
=:Bb+Bs+Sb+Ss\ .\qquad\qquad\qquad\qquad\qquad\qquad\qquad\quad\cr
}
$$
Let us now bound separately each of the four terms $Bb,\; Ss,\; Sb,\; Bs$.
%
Using (\ref{mom})
\begin{align*}
Bb&\le\sum_{j=0}^{{J}}2^{j(\gamma-d)}\bE\sum_{\eta\in \Z_j}  |\hat\bjk-\bjk|^z
1_{\{|\bjk|\ge \frac\kappa 2  t_n{}\}}\le\\
&\le
\sum_{j=0}^{{J}}\sum_{\eta\in \Z_j}1_{\{ |\bjk|\ge \frac\kappa 2    t_n{}\}}
2^{j(\gamma-d)}\bE  |\hat\bjk-\bjk|^z \le\\
&\le
\sum_{j=0}^{{J}}\sum_{\eta\in \Z_j}1_{\{ |\bjk|\ge \frac\kappa 2    t_n{}\}}
2^{j(\gamma-d)}s_zn^{-\frac z2}\le\quad (*)
\\
&\le C2^{J_1\gamma}n^{-\frac z2}
\end{align*}
where $J_1$ again is chosen such that
for $j\ge J_1$, $|\bjk|\le \frac\kappa 2 t_n$. To prove (\ref{propsum1}), we stop in (*), the next bound yields (\ref{propsum}).
\begin{align*}
Ss&\le\sum_{j=0}^{{J}}2^{j(\gamma-d)}\sum_{\eta\in \Z_j}  |\bjk|^z
1_{\{|\bjk|< 2{\kappa}  t_n{}\}}\quad (*)
\\
&\le\sum_{j=0}^{{J_1}}2^{j(\gamma)}[2{\kappa}  t_n]^z
+\sum_{j=J_1+1}^{{J}}2^{j(\gamma-d)}\sum_{\eta\in \Z_j}  |\bjk|^z
%
\end{align*}
which gives the proper rate of convergence. Again, to prove (\ref{propsum1}), we stop in (*), the next bound yields (\ref{propsum}).
Moreover, using (\ref{momsup}) and (\ref{gd}),
\begin{align*}
Bs&\le \sum_{j=0}^{{J}}2^{j(\gamma-d)}\bE\sum_{\eta\in \Z_j}  |\hat\bjk-\bjk|^z
1_{\{|\hat\bjk-\bjk|\ge \frac\kappa 2   t_n{}\}}
1_{\{|\bjk|< \frac\kappa 2  t_n{}\}}\\
&\le \sum_{j=0}^{{J}}2^{j(\gamma-d)}\bE\sum_{\eta\in \Z_j}  |\hat\bjk-\bjk|^z1_{\{|\hat\bjk-\bjk|\ge \frac\kappa 2   t_n{}\}}\\
&\le \sum_{j=0}^{{J}}2^{j(\gamma-d)}\sum_{\eta\in \Z_j} [\bE |\hat\bjk-\bjk|^{2z}]^{\frac12}\bP\{|\hat\bjk-\bjk|\ge \frac\kappa 2   t_n{}\}^{\frac12}\\
&\le \sum_{j=0}^{{J}}2^{j(\gamma)} [s_{2z}n^{-z}]^{\frac12}[c2^{jd}n^{-6\kappa}]^{\frac12}\le Cn^{-\frac z2}
\end{align*}
for $\kappa>\frac \gamma{6d}$.
Finally, using
(\ref{gd}), and the fact that for $f$ bounded, $|\bjk|\le C2^{-\frac{jd}2}$
\begin{align*}
Sb&\le \sum_{j=0}^{{J}}2^{j(\gamma-d)}\bE\sum_{\eta\in \Z_j}  |\bjk|^z
1_{\{|\bjk-\hat\bjk|\ge {\kappa}  t_n{}\}}1_{\{|\bjk|\ge 2{\kappa}  t_n{}\}}
\\
&\le \sum_{j=0}^{{J}}2^{j(\gamma-d)}M2^{-jz\frac d2}\bP\{|\bjk-\hat\bjk|\ge {\kappa}  t_n{}\}
\\
&\le \sum_{j=0}^{{J}}[c2^{j[-\frac {dz}2+\gamma]}n^{-6\kappa}]
\\
&\le C [2^{J[-\frac {dz}2+\gamma]}n^{-6\kappa}]\le Cn^{\frac 12}
,\end{align*}
for $\kappa>0$.
%
\section{Proof of the lower bound}
Let us recall that given
two probabilities $P$, $Q$ on some measure space their Kullback-Leibler distance
is
$$
K(P,Q)=\left \{
 \begin{array}{ll}
   \int \log \frac{dP}{dQ} dP =\int \frac{dP}{dQ}\,\log \frac{dP}{dQ}\,  dQ & \mbox{ if } P<<Q \\
  \\
 +\infty &  \mbox{otherwise\ .}
  \end{array}
\right.
$$
If $P$, $Q$ are probabilities on $\mathbb{S}^d$ having densities $f$, $g$ respectively
with respect to
Lebesgue measure, then if $g$ is bounded below by some constant $c>0$
\begin{equation}\label{bound-kullback-l2}
\begin{array}{c}
\displaystyle K(P,Q)=\int \log\tfrac fg\, f\, dx=\int \log(\tfrac {f-g}g+1)\, f\, dx
\le\cr\le \int \tfrac {f-g}g\, f\, dx
\displaystyle =\int \tfrac {(f-g)^2}g\, dx\le \tfrac 1c \Vert f-g\Vert_2^2\ .
\end{array}
\end{equation}
We make use of Fano's lemma below, see \cite{tsyb}  and the references therein.
We use the point of view  introduced in  \cite{birg}.
\begin{theorem} (Fano's Lemma)
Let  $\mathscr{F}$ be a sigma algebra on the space  $\Omega.$
Let $F_i \in \mathscr{F}, ~~ i\in\{ 0,1,\ldots,m \}$ such that $\forall i\neq j,
F_i\cap F_j =\emptyset.$
Let $P_i$, $i=0,\ldots,m$ be probability measures on $(\Omega,\mathscr{A}).$
If
\begin{align}
&p \enspace\mathop{=}^{def}\enspace\sup _{i=0,\dots, m} P_{i}(  F^c_i),\\
&\kappa(P_0,\dots,P_m) \enspace\mathop{=}^{def}\enspace\inf_{j=0,\dots,m}
\frac 1m \sum_{i\neq j} K(P_i,P_j)\ .
\end{align}
then
 \begin{equation}
p \geq \frac 12 \wedge C (\sqrt m \,    e^{ - \kappa(P_0,\dots,P_m)}),\quad C= e^{3/e}\ .
\end{equation}
\end{theorem}
$\bullet$ We prove first that the minimax $L^p$-loss is $\gtrsim n^{-\alpha p}$ with
$\alpha=\frac s{2s+d}$.
For every $j$ let us consider the family $\ccl A_j$ of densities
\begin{equation}\label{densities}
f_\varepsilon=\frac 1{|\mathbb{S}^d|}+\gamma \sum_{\xi\in A_j}\varepsilon_\xi\psi_{j,\xi}
\end{equation}
where $A_j$ is a subset of $\ccl Z_j$ to be made precise later,
$\varepsilon_\xi\in\{0,1\}$ and $\gamma$ is chosen so that all these
functions are positive. We are going to show that for every estimator $\hat f$,
$$
\sup_{f_\varepsilon}\E_{f_\varepsilon}\Vert\hat f-f_\varepsilon\Vert_p^p\ge c n^{-\frac {sp}{2s+d}}\ .
$$
Throughout this section we shall note $x\lesssim y$, $x\gtrsim y$
whenever it holds $x\le c y$ or $x\ge c y$ respectively, $c$ being a
strictly positive constant independent of $j,\xi$. We shall note
$x\simeq y$ whenever both $x\lesssim y$ and $x\gtrsim y$ hold.
Thanks to (\ref{lp}) for these functions to be positive it
is enough that $|\gamma|\lesssim 2^{-jd/2}$. Such a $\gamma$ can
even be chosen in such a way that all the densities (\ref{densities}) are bounded from below
by a strictly positive constant. If the functions $(\psi_{j\xi})_{\xi\in\ccl Z_j}$ were
orthonormal we would have immediately that
$$
\Big\Vert\sum_{\xi\in \ccl Z_j}\lambda_\xi\psi_{j,\xi}\Big\Vert_p\ge
c\Bigl(\sum_{\xi\in \ccl Z_j}|\lambda_\xi|^p\Vert\psi_{j,\xi}\Vert_p^p\Bigr)^{1/p}
$$
Needlets are not a basis, but their scalar product is close enough to $0$ if the respective
cubature points are far enough. Hence one can get the following lemma that
states that a subset $A_j\subset \ccl Z_j$
can be chosen so that it is quite large and inequalities (\ref{inequality2}) and (\ref{inequality}),
in a sense, can be reversed.
\begin{lemma}\label{reverse-inequality}
There exists a subset $A_j\subset \ccl Z_j$ such that ${\rm card}\, A_j\gtrsim 2^{jd}$ and
$$
\Big\Vert\sum_{\xi\in A_j}\lambda_\xi\psi_{j,\xi}\Big\Vert_p\ge
\begin{cases}
c\sup_{\xi\in A_j}|\lambda_\xi|\Vert\psi_{j,\xi}\Vert_\infty&\mbox{if }p=\infty\cr
c\Bigl(\sum_{\xi\in A_j}|\lambda_\xi|^p\Vert\psi_{j,\xi}\Vert_p^p\Bigr)^{1/p}
&\mbox{if }p<+\infty\ .
\end{cases}
$$
\end{lemma}
Let us now impose conditions that ensure that $f_\varepsilon$ belongs to the ball
$B^s_{r,q}(M)$. Now, recalling (\ref{bnorm}),
$$
\Vert f_\ep\Vert_{B^s_{r,q}}=|\gamma|2^{j(s+d(\frac12-\frac1r))}
\Bigl(\sum_{\xi\in\ccl Z_j}|\ep_\xi|^r\Bigr)^{1/r}\lesssim
|\gamma|2^{j(s+d(\frac12-\frac1r))}2^{jd/r}
$$
where we use the fact that $|\varepsilon_\xi|=1$. Therefore the condition
$\Vert f_\ep\Vert_{B^s_{r,q}}\le M$ follows from
$$
|\gamma|\lesssim M2^{-j(s+\frac d2)}\ .
$$
In order to apply Fano's Lemma and get a lower bound of the
left hand side let us first get an upper bound for the Kullback-Leibler distances
$K(f_\varepsilon;f_\varepsilon)$, which comes from (\ref{bound-kullback-l2}) and
(\ref{lp}) for $p=2$,
\begin{equation}\label{lb-bornel2}
\Vert f_\varepsilon-f_{\varepsilon'}\Vert_2^2\le \gamma^2\sum_{\xi\in A_j}
|\varepsilon_\xi-\varepsilon'_\xi|^2<\gamma^22^{jd}\le 2^{-2js}\ .
\end{equation}
%
By Lemma \ref{reverse-inequality}
$$
\Vert f_\varepsilon-f_{\varepsilon'}\Vert_p
\ge \Bigl(\sum_{\xi\in A_j}|\varepsilon_\xi-\varepsilon'_\xi|^p
\Vert\psi_{j,\xi}\Vert_p^p\Bigr)^{1/p}\ .
$$
Thanks to Lemma \ref{reverse-inequality}, the set of functions $\ccl A_j$
has a cardinality that is $\ge 2^{c2^{jd}}$. By the Varshanov-Gilbert Lemma
(\cite{tsyb} e.g.) there exists a subset $\ccl A'_j\subset \ccl A_j$
such that ${\rm card}\, \ccl A'_j\ge 2^{c'2^{jd}}$ and such that if
$f_\varepsilon,f_\varepsilon'\in \ccl A'_j$, then
$\sum_{\xi\in A'_j}|\varepsilon_\xi-\varepsilon'_\xi|>\frac 14\,2^{jd}$. Therefore,
as $|\varepsilon_\xi-\varepsilon'_\xi|$ can be $=0$ or $=1$ only and by (\ref{lp}),
$$
\Vert f_\varepsilon-f_{\varepsilon'}\Vert_p\gtrsim |\gamma|
2^{jd(\frac 12-\frac 1p)}(\tfrac {2^{jd}}4)^{1/p}\simeq 2^{-js}
$$
which implies that the events
$\{\Vert \widehat f-f_{\varepsilon'}\Vert_p\ge \tfrac 12\,2^{-js}\}$ are disjoint.
The family of densities $f_\varepsilon\in \ccl A'_j$ given by the
Varshanov-Gilbert Lemma has cardinality $m\simeq2^{c'2^{jd}}$ and by (\ref{lb-kbl2})
and (\ref{lb-bornel2})
$$
K(f_\varepsilon,f_{\varepsilon'})\lesssim
\Vert f_\varepsilon-f_{\varepsilon'}\Vert_2^2\lesssim
2^{-2js}\ .
$$
We apply now Fano's lemma to the probabilities $\P_\ep$ that are the
$n$ times product of $f_\ep\, dx$ and to the events
$A_\ep=\{\Vert \widehat f-f_{\varepsilon'}\Vert_p\ge \tfrac 12\,2^{-js}\}$. It is well known that
$$
K(\P_\ep;\P_{\ep'})=nK(f_\ep;f_{\ep'})\ .
$$
By Markov inequality and Fano's lemma
$$
\displaylines{
\sup_{f_\varepsilon\in\ccl A_j'}
\E\Vert \hat f-f_\varepsilon\Vert_p^p\ge 2^{-p}2^{-jsp}\sup_{i\le m}
\P_\ep(\Vert \hat f-f_\varepsilon\Vert_p>\delta)\gtrsim\cr
\gtrsim 2^{-jsp}
\bigl(\tfrac 12\,\wedge e^{-n2^{-2js}}
\underbrace{\sqrt{\#\ccl A_j}}_{\simeq 2^{c2^{jd}}}\,\bigr)\ .\cr
}
$$
Now let  $j$ be so that
$
n2^{-2js}\simeq 2^{jd}
$,
that is $2^j\simeq n^{\frac 1{2s+d}}$. With this choice one has
$$
\tfrac 12\,\wedge (e^{-n2^{-2js}} e^{c2^{jd}})\ge c>0\ .
$$
Therefore
$$
\sup_{f_\varepsilon\in\ccl A_j'}
\E\Vert \hat f-f_\varepsilon\Vert_p^p\gtrsim c2^{-jsp}\sim
n^{-\frac {sp}{2s+d}}\ .
$$
$\bullet$ We prove now that the minimax $L^p$-loss is $\gtrsim n^{-\frac {p(s+d(\frac 1p-\frac 1r))}{2(s+d(\frac 12-\frac 1r))}}$.
Let us consider the two densities
$$
f_0=\frac 1{|\mathbb{S}^d|}+\gamma \psi_{j,\xi},\quad f_0=\frac 1{|\mathbb{S}^d|}+\gamma \psi_{j,\xi'}
$$
with $\gamma$ such that the above are positive ($|\gamma|\lesssim 2^{-jd/2}$ is enough). If
$|\gamma|\le 2^{-js} 2^{-jd(\frac 12-\frac 1r)}M$, then thanks to (\ref{bnorm})
both $f_0$ and $f_1$ belong to the ball $B^s_{r,q}(M)$. Remark that this condition implies
$|\gamma|\lesssim 2^{-jd/2}$, as we assume $s\ge \frac dr$.
Also
$$
K(f_0\, dx,f_1\, dx)\le \Vert f_0-f_1\Vert_2^2\approx \gamma^2
$$
so that,if we denote by $\P_0,P_1$ the $n$-times product of $f_0\, dx$ and $f_1\, dx$
by itself respectively,
$K(\P_0,\P_1)\approx n\gamma^2$. By \tref{inequality} and Lemma \ref {reverse-inequality} we have
\begin{equation}\label{lb-kbl2}
\begin{array}{c}
\displaystyle\Vert f_0-f_1\Vert_p=|\gamma|\Vert \psi_{j,\xi}-\psi_{j,\xi'}\Vert_p\ge |\gamma|
2^{jd(\frac 12-\frac 1p)}\sim\\
\displaystyle\vphantom{\sum^2}\sim 2^{-j(s+d(\frac 12-\frac 1r)}2^{jd(\frac 12-\frac 1p)}=r
=2^{-j(s+d(\frac 1p-\frac 1r))}\ .\\
\end{array}
\end{equation}
We choose $\gamma=\frac 1{\sqrt{n}}=2^{-j(s+d(\frac 12-\frac 1r))}$, so that
$K(\P_0,\P_1)\approx n$. Moreover with this choice of $n$, $j\approx \log n((2(s+d(\frac 12-\frac 1r)))^{-1}$,
so that again by Fano's lemma,
$$
\sup_{i=1,2}\E\Vert \hat f-f_i\Vert_p^p\ge \delta^p\sup_{i=1,2}
\P_i(\Vert \hat f-f_i\Vert_p\ge \delta)\ .
$$
Thanks to \tref{lb-kbl2} the events $\{\Vert \hat f-f_i\Vert_p\ge \delta\}$
are disjoint if $\delta\lesssim 2^{-j(s+d(\frac 1p-\frac 1r))}$. Therefore by
Fano's lemma
$$
\sup_{i=1,2}\E\Vert \hat f-f_i\Vert_p^p\gtrsim \delta^p \gtrsim
2^{-j(s+d(\frac 1p-\frac 1r))}=
n^{-\frac {p(s+d(\frac 1p-\frac 1r))}{2(s+d(\frac 12-\frac 1r))}}\ .
$$
We have therefore proved that
$\sup_{f\in B^s_{rq}(M)}\min_{\hat f}\EE_f(\|\hat f-f\|_p^p$ is
$$
\gtrsim n^{-\frac {sp}{2s+d}}  \quad\mbox{and}\quad \gtrsim
n^{-\frac {p(s+d(\frac 1p-\frac 1r))}{2(s+d(\frac 12-\frac 1r))}}\ .
$$
Putting things together and checking for which values of the parameters one rate is
larger than the other one concludes the proof of Theorem \ref{lb-fin}.
Note that, as $s>\frac dr$, if $1\le p\le 2$
$$
\frac {sp}{2s+d}\le \frac {p(s+d(\frac 1p-\frac 1r))}{2(s+d(\frac 12-\frac 1r))}\cdotp
$$
%

%
%

\bibliography{pencho}
\bibliographystyle{alpha}

{\parindent0pt
\hbox{
\vtop{\hsize 6cm\footnotesize
PB \& DM

Dipartimento di Matematica

Universit\`a di Roma {\it Tor Vergata}

Via della Ricerca Scientifica

00161 Roma (Italy)
\medskip



baldi@mat.uniroma2.it

marinucc@mat.uniroma2.it
}
\qquad
\vtop{\hsize 7cm
\footnotesize
GK \& DP

Laboratoire de Probabilit\'es et
Mod\`eles Al\'eatoires

2, Pl. Jussieu

75251 Paris Cedex 05, France
\medskip



kerk@math.jussieu.fr,

picard@math.jussieu.fr
}
}
}

\end{document}